%% file: main1.tex
\documentclass[a4paper,12pt]{article}
\usepackage[T2A]{fontenc}
\usepackage[cp1251]{inputenc}
\usepackage{amsmath}
\usepackage{amssymb}
\usepackage[dvips]{graphicx}
\usepackage[english,russian]{babel}
\usepackage{ifpdf}

\newtheorem{theorem}{Теорема}[section]
\newtheorem{corollary}[theorem]{Следствие}
\newtheorem{lemma}[theorem]{Лемма}
\newtheorem{proposition}[theorem]{Предложение}

\newtheorem{definitionhead}[theorem]{Определение}
\newenvironment{definition}{\begin{definitionhead}%
\sl}{\end{definitionhead}}

\newcommand\Proof{\noindent{\bf Доказательство. }}
\newcommand\Endproof{\nopagebreak\strut%
\nopagebreak\hfill\nopagebreak$\Box$\medbreak}
\newcommand\Noproof{\nopagebreak\strut%
\nopagebreak\hfill\nopagebreak$\Box$\par}

\def\Pref{\operatorname{Pref}}
\def\Suf{\operatorname{Suf}}
\def\PAL{\operatorname{PAL}}
\def\Fol{\operatorname{Fol}}
\def\Card{\operatorname{Card}}
\def\Red{\operatorname{Red}}
\def\Ret{\operatorname{Ret}}
\def\mod{\operatorname{mod}}
\def\dim{\operatorname{dim}}
\def\Const{\operatorname{Const}}

\title{Нормальные базисы и символическая динамика}
\author{А.~Л.~Чернятьев}
\date{}
\begin{document}
\fontsize{14pt}{18} \selectfont

\pagestyle{empty}
\begin{center}
\vfill\vfill {\LARGE Московский государственный университет

им. М.В.~Ломоносова}

\vfill \hfill На правах рукописи

\hfill УДК 512+519.17

\vfill\vfill {\LARGE Чернятьев Александр Леонидович}

\vfill{\huge Нормальные базисы и символическая динамика}

\vfill Специальность 01.01.06~--- математическая логика,

алгебра и теория чисел

\vfill Диссертация на соискание ученой степени кандидата
физико-математических наук

\vfill Научные руководители:

доктор физико-математических наук,

профессор А.В. Михалев,

доктор физико-математических наук,

А.Я. Белов

\vfill\vfill\vfill Москва

 2008
\end{center}
\newpage
%\fontsize{14pt}{17}
\pagestyle{plain} \tableofcontents
\newpage

{\bf }

\input intro.tex

\newpage

\input chapter1.tex

\newpage

\input chapter2.tex

\newpage

\input chapter3.tex

\newpage

\input chapter4.tex

\newpage

\input bib.tex
\end{document}

%% file: intro.tex
\section{Введение.}
\paragraph{Актуальность темы.}
Комбинаторика слов находит свое применение в самых разных разделах
математики. Например, в алгебре при изучении базисов и нормальных
форм, в алгебраической топологии, в символической динамике. Ряд
проблем, относящихся к комбинаторике слов находится на стыке алгебры
и теории динамических систем. Многие проблемы комбинаторики слов
представляют самостоятельный интерес.

Комбинаторика слов широко используется  в задачах комбинаторной
теории групп (\cite{LSh},\cite{NA}), в теории алгебр Ли, вопросах
бернсайдовского типа  и задачах, связанных с мономиальными
алгебрами. Комбинаторная техника, относящаяся к теории групп,
развивалась в работах М.~Дэна, Е.~С.~ Голода и И.~Р.~ Шафаревича,
П.~С.~Новикова, С.~И.~ Адяна, И.~Рипса, М.~Громова,
А.~Ю.~Ольшанского, М.~В.~Сапира.

Е.~С.~ Голод и И.~Р.~ Шафаревич \cite{GolShaf}, \cite{Golod2}
построили конечно порожденную бесконечную периодическую группу (с
неограниченной экспонентой) на основе рассмотрения нормальных форм
алгебр и оценки функий роста. П.~С.~Новиков и С.~И.~Адян \cite{NA}
провели детальное исследование свойств периодичности, находящее свое
применение в самых разных разделах математики. Ими были впервые
построены примеры бесконечных конечно порожденных периодических
групп ограниченой экспоненты (т.е. решена проблема Бернсайда),
получены наилучшие из известных оценок на экспоненту для таких групп
(см. обзор \cite{AdianRasb}). В дальнейшем был исследован случай
четной экспоненты (см. \cite{Lys}, \cite{Ivanov}).

В основе замечательных работ М.~Громова и А.~Ю.~Ольшанского также
лежит техника диаграмм Ван-Кампена, возникшая в комбинаторной
топологии. Подробнее (а также литературу по этой теме) см. в
монографии \cite{OlshanskiBook}.

Комбинаторные соображения, возникшие в символической динамике
(автоматные группы), нашли свое применение в работах С.~В.~Алешина
\cite{Alesh1,Alesh2,Alesh3} и и Р.~И.~Григорчука \cite{Grig},
\cite{Grig1}, \cite{Grig2} при решении проблемы Милнора --
посторении групп промежуточного роста (при этом группы Григорчука
периодичны). Впервые автоматные полугруппы были построены в работах
С.~В.~Алешина.  (Изложение примера С.~В.~Алешина -- см. в книге
\cite{KargMerzl}) Автоматные конструкции активно используются в
самых разных ситуациях (см. \cite{BBL}, \cite{Uf}, \cite{Salomaa},
\cite{KudriavAleshPodkolzin}). Возникают они и у нас (графы Рози).

Комбинаторика слов активно используется в алгебрах Ли, особенно в
проблемах бернсайдовского типа \cite{Kostrikin}. В теории алгебр Ли
описание базиса дается в терминах так называемых ``правильных слов''
(базис Холла-Ширшова). Слово называется {\it правильным} если оно
лексикографически больше любого его циклически сопряженного (слова
$v_1,v_2$ {\it циклически сопряжены}, если $v_1=u_1u_2$,
$v_2=u_2u_1$ для некоторых $u_1$ и $u_2$). (А запись слова по циклу
используется в теории групп, тесно связанной с теорией алгебр Ли.)
Только в правильном слове (причем однозначно) можно расставить лиевы
скобки так, чтобы при их раскрытии исходное слово оказалось старшим
членом получившегося полинома. Тем самым правильные слова задают
базис свободной алгебры Ли (базис Холла--Ширшова) \cite{Bach}.
Применив методы символической динамики (равномерно рекуррентные
слова и соображения компактности) Д.~Бэкелин установил, что любое
сверхслово содержит подслово вида $uvu$, где $u$ и $v$~-- правильные
слова, получив, тем самым, короткое доказательство локальной
нильпотентности подалгебры алгебры Ли, порожденной сэндвичами,
упростив соответствующие работы А.~И.~Кострикина \cite{Uf}.

Применив гомологическое соображение, связанное с невозможностью
зацепления правильного слова за самого себя, А.~И.~Ширшов показал
алгоритмическую разрешимость проблемы равенства в алгебрах Ли с
одним определяющим соотношением.

С помощью комбинаторики слов А.~И.~Ширшов \cite{Shirshov4} доказал
теорему о свободе подалгебры свободной алгебры Ли. Для супералгебр
это обобщил А.~А.~Михалев \cite{MikhSmall1}, \cite{MikhSmall3}. Он
показал алгоритмическую разрешимость проблемы равенства для алгебр
ли с одним определяющим соотношением.

Комбинаторика слов активно используется в проблемах бернсайдовского
типа, в теории $PI$-алгебр, достаточно упомянуть знаменитую теорему
Ширшова о высоте \cite{Shirshov1}, \cite{Shirshov2} утверждающую
возможность приведения слов к кусочно периодическому виду.

{\bf Теорема А.И.Ширшова о высоте.} {\it Пусть $A$ -- конечно
порожденная $PI$-алгебра степени $m$. Тогда существует конечный
набор элементов $Y$ и число $h\in \mathbb{N}$ такие, что $A$ линейно
представима (то есть порождается линейными комбинациями) множеством
элементов вида:
\begin{center}
$w={u_1}^{k_1}{u_2}^{k_2}\cdots {u_r}^{k_r}$,  где $u_i\in Y$ и $r
\leq h$.
\end{center}
}

При этом в основе оригинальных доказательств А.~И.~Ширшова  (как
теоремы о свободе так и теоремы о высоте) лежала техника, связанная
с преобразованием алфавита путем подстановок. Эта же техника
используется при работе с равномерно--рекуррентными словами и в
символической динамике.

Последующие доказательства теоремы о высоте \cite{Belov} и ее
обобщение для алгебр Ли \cite{Mishenko1} использовали анализ свойств
периодичности.

Понятие {\it роста} в алгебре является важным комбинаторным
инвариантом, ему посвящена монография \cite{KL} (см. также
\cite{Uf}, \cite{BBL}).
 Если размерность пространства, порожденного словами степени не выше $n$ от
 образующих $A$ растет как $n^\lambda$, то величина $\lambda$ называется
 {\it размерностью Гельфанда--Кириллова алгебры $A$}. Размерность
  Гельфанда--Кириллова может быть равной $0$, $1$, а также
  любому числу $\ge 2$, $\infty$ или не существовать.
  То обстоятельство, что она не может принимать промежуточные
  значения на интервале $(1,2)$  составляет содержание знаменитой
  {\it Bergman gap theorem.} Ассоциативная алгебра размерности Гельфанда--Кириллова
  0 конечномерна. Л.~Смолл показал, что ассоциативная алгебра
  размерности Гельфанда--Кириллова 1 является $PI$-алгеброй.
  Базисы ассоциативных алгебр размерности Гельфанда--Кириллова
  больше 1 с минимальной функцией роста исследовались в
  работах А.~Т.~Колотова \cite{Kol}, \cite{Kol1}. Их описание
  дается в терминах так называемых {\it последовательностей Штурма}
  находящейся в центре внимания данной работы.

Обобщение понятия роста на бесконечномерный случай является понятие
{\it ряда коразмерностей}, введенное А.~Регевым. Первоначальное
доказательство А.~Регева об экспоненциальной оценке ряда
коразмерности относительно свободных алгебр было улучшено
В.~Н.~Латышевым с помощью оценки числа полилинейных $n$-разбиваемых
слов на основе теоремы Дилуорса \cite{Latyshev12}. Само же понятие
{\it $n$-разбиваемого слова}
 возникло у А.~И.~Ширшова в его теореме о высоте.
Ряды коразмерности исследовались также в работах В.~Н.~Латышева,
С.~П.~Мищенко, М.~В.~Зайцева, A.~Giambruno. Понятию роста посвящена
монография \cite{KL}.

Комбинаторика слов работает в теории полугрупп. Следует указать
работы Екатеринбургской школы Л.~Н.~Шеврина, в часности, работы
М.~В.~ Сапира, О.~Г.~ Харлампович. Они активно применяли методы
символической динамике к теории полугрупп.

В теории мономиальных алгебр комбинаторика слов имеет
основополагающее значение и находится в центре внимания работы
\cite{BBL}.

Структурная теория позволила получить элегантные, но
неконструктивные доказательства в теории колец. Вместе с тем она
оказала тормозящее влияние на развитие непосредственно комбинаторных
методов, пусть более трудоемких, но зато позволяющих получать
конструктивные оценки и дающих лучшее понимание комбинаторной сути.

Вопросы, связанные с базисами алгебр, приводят изучению бесконечных
(в одну или обе стороны) слов или {\it сверхслов}. Фундаментальным
понятием в теории сверхслов является понятие {\it
равномерно-рекуррентного} слова, введенное Х.~Фюрстенбергом
\cite{F}. Слово $W$ называется {\it равномерно-реккурентным}, если
для каждого подслова $v\subset W$ существует натуральное $N(v)$,
такое, что для любого подслова $u\subset W$  длины не менее, чем
$N(v)$, $v$ является подсловом $u$.

Имеет место следующая \par {\bf Теорема}. {\it Пусть $W$ --
бесконечное сверхслово. Тогда существует такое
равномерно-рекуррентное слово $\hat{W}$, все подслова которого
являются подсловами $W$.}

Эта теорема исключительно важна в комбинаторике слов, поскольку
очень часто позволяет свести изучение произвольных слов к изучению
равномерно-реккурентных слов.

В терминах равномерно-рекуррентных слов строится теория радикалов
мономиальных алгебр (\cite{BBL}). Мономиальная алгебра называется
{\it мономиально почти простой}, если фактор по идеалу, порожденному
по любому моному, нильпотентен.

Множество ненулевых слов в почти простой мономиальной алгебре есть
множество всех подслов некоторого равномероно-рекуррентного слова.

Пересечение же идеалов с почти простым фактором совпадает с
нильрадикалом мономиальной алгебры, а также ее радикалом Джекобсона
(см. \cite{BG}).

В терминах равномерно-рекуррентных слов также получается описание
слабо нетеровых мономиальных алгебр. Каждое ненулевое слово слабо
нетеровой мономиальной алгебры есть подслово из набора (сверх)слов,
удовлетворяющего следующему условию: каждое слово из этого набора
либо конечное, либо является бесконечным (односторонними или
двухсторонними) словом, которое при выбрасывании некоторого
конечного куска распадается на равномерно-рекуррентные части.

Существует разные подходы к изучению сверхслов:

\begin{enumerate}
\item Непосредственно комбинаторные свойства слов

\item Графы подслов, или {\it графы Рози}

\item Топологическая динамика
\end{enumerate}

%Библиография по общей теории слов.
Классическими работами по теории комбинаторики слов являются
монографии \cite{L1,L2}, \cite{RS}.

Другим инструментом изучения сверхслов является понятие графов
подслов или {\it графов Рози}. Если $W$ -- бесконечное сверхслово
над алфавитом $A$, то {\it $k$-графом Рози} называется граф, вершины
которого соответствуют различным подсловам $W$ длины $k$. Из вершины
$w_1$ в вершину $w_2$ ведет стрелка, если максимальный суффикс $w_1$
совпадает с максимальным префиксом $w_2$, то есть $w_1=a_1u$,
$w_2=ua_2$, где $a_1,a_2 \in A$.

Общий подход, связанный с описанием слов с помощью динамических
систем, следующий.  Пусть $W=\{w_n\}$ -- бесконечное слово. $\tau
(\{w_n\})=\{w_{n+1}\}$ -- оператор сдвига. Рассмотрим замыкание
траектории слова относительно метрики Хэмминга $X \in A^*$. Прямые
задачи символической динамики связаны с получением информации о
динамической системе $(X,\tau)$ по информации о слове $W$.

Известно, что если слово $W$ {\it равномерно-рекуррентно}, то
полученная динамическая система минимальна, то есть не содержит
нетривиальных замкнутых инвариантных траекторий.

Также стоит отметить работу  \cite{BK}, изучающую слова, получаемые
из взятия дробных частей многочленов со старшим иррациональным
коэффициентом в целых точках.

Обратно, пусть имеется дискретная топологическая динамическая
система, то есть задано компактное топологическое пространство $M$,
гомеоморфизм $f:M\to M$ и несколько открытых подмножеств
$$
U_1,U_2,\ldots, U_{n-1}.
$$
Положим также
$$
U_n=M \setminus U_1\cup U_2\cup \ldots \cup U_{n-1}.
$$
Рассмотрим начальную точку $x \in M$ и последовательность итераций
$f^{(-1)}(x),x,f(x),\ldots$. По этой последовательности можно
построить слово $W=\{w_n\}$ над алфавитом $A=\{a_1,a_2,\ldots,
a_n\}$ следующим образом: \par $w_i=a_k$, если $f^{(i)}(x) \in U_k$.
По свойствам динамической системы (размерность множества $M$,
эргодичность) можно получить информацию о слове $W$.

Важным примером в изучении сверхслов динамическим подходом являются
слова с предельной функцией роста.

Хорошо известно, что если функция роста слова $V(n)$ (то есть
размерность пространства, порожденного словами степени не выше $n$)
при некотором $n$ удовлетворяет неравенству $V(n)<n(n+3)/2$, то
алгебра имеет линейный рост.

В работе Колотова \cite{Kol} \label{Kol} построены алгебры с
``предельной'' функцией роста (то есть когда $V(n)=n(n+3)/2$),
которые описаны в терминах поворота окружности. А именно, все такие
алгебры, кроме счетного множества, строятся как алгебры $A_W$, где
$W=\{w_i\}$ -- слово над алфавитом $\{0,1\}$, задаваемая
иррациональными $\alpha,\beta\in (0,1))$: $w_i=g(i+1)-g(i)$, где
$g(i)=[\alpha i + \beta]$. В комбинаторике слов чаще используется
{\it функция сложности}, $T(n)$ равная количеству различных подслов
длины $n$. И, таким образом, $V(n)=\sum_{k}T_k$. Для алгебр функцию
сложности корректно будет определить следующим образом:
$T_A{n}=V_A(n)-V_A(n-1)$, поскольку алгебра может быть неоднородна.

Известно, что $\lim_{n\to\infty} T_A(n)-n$ всегда существует. Он
может быть равен $-\infty, C, +\infty$.   Если $\lim_{n\to\infty}
T_A(n)-n=-\infty$, то алгебра $A$ либо конечномерна, либо имеет
медленный рост. Л.Смолл и Д.Бергман исследовали алгебры медленного
роста в ряде своих работ. Суммируя их результаты, получается
описание нормальных базисов таких алгебр.

Назовем алгебру {\it граничной}, если $\lim_{n\to\infty}
T_A(n)-n=C$. Описание нормальных базисов граничных алгебр следует из
теоремы \ref{MinGrow} и результатов \cite{BBL}.

Слова с предельной функцией роста $T(n)=n+1$ образуют класс так
называемых {\it слов Штурма (Sturmian words)}, другое название --
{\it слова Бетти (Beaty words)}, которые впервые были упомянуты в
работе \label{HedlundMorse} \cite{MorseHedlund} Hedlund, Morse
``Symbolic dynamics II. Sturmian trajectories''.

Классическая теория слов Штурма описана в обзорах
\label{BerstelSeebold} \cite{BS}, \cite{AdL}.

К наиболее важным результатам в теории слов Штурма относится так
называемая {\it теорема эквивалентности}, в которой утверждается
эквивалентность трех классов сверхслов над двубуквенным
алфавитом:\par {\bf Теорема
эквивалентности.}(\cite{MorseHedlund})\label{TheorEq} {\it Следующие
условия на слово $W$ эквивалентны:
\begin{enumerate}
\item Слово $W$ имеет фунцию сложности $T_W(n)=n+1$.

\item Слово $W$ сбалансированно и непериодично.

\item Слово $W$ порождается системой $(\mathbb{S}^1, U, T_\alpha)$
с иррациональным углом вращения $\alpha$.
\end{enumerate}
}

Последние продвижения в теории слов Штурма описаны в обзоре
\label{BerstelRecent} \cite{B} J. Berstel ``Recent results in
Sturmian words''.

Естественными обобщениями слов Штурма являются слова с минимальным
ростом, то есть слова с функцией роста $T(n)=n+K$, начиная с
некоторого $n$. Для двубуквенных алфавитов они носят название {\it
квазиштурмовых} слов. Слова с функцией роста, удовлетворяющей
соотношению $\lim_{n \to \infty} T(n)/n = 1$ изучены в работе
\label{Aberkane} \cite{Ab} A.Aberkane ``Words whose complexity
satisfies $\lim p(n)/n = 1$''.

Другим обобщением слов Штурма является обобщение, связанное с
понятием {\it сбалансированности}, а также {\it
$m$-сбалансированности}. Сбалансированные непериодические слова над
$n$-буквенным алфавитом изучены в работе \label{Graham} \cite{GR}
Graham ``Covering the Positive Integers by disjoints sets of the
form $\{n\alpha + \beta\}:n=1; 2;\ldots$''. Построение порождающей
динамической системы для сбалансированных непериодических слов
является одним из результатов данной работы. Исследование сбалансированных слов тесно
связано с построением ненильпотнентных ниль-алгебр (см. $3.9$)

Описание периодических сбалансированных слов связано с {\it
гипотезой Френкеля (Fraenkel's conjecture)}, утверждающей, что все
сбалансированные периодические слова над алфавитом $A=\{a_1,\ldots,
a_n\}$ из $n$ символов имеют вид
$$
W=(U_n)^{\infty},
$$
где $U_n$ задается рекуррентно:
$$
U_n=(U_{n-1}a_nU_{n-1}), \ U_3=a_1a_2a_1a_3a_1a_2a_1.\\
$$
Для 3-х буквенного алфавита гипотеза была доказана Р.~Тайдеманом
(\cite{T1,T2}). В настоящий момент гипотеза доказана для алфавитов
состоящих не более чем из $7$ символов.

Продвижение в задачах символической динамики для слов с линейной
функцией роста получено в работе \cite{AR} \label{ArnouxRazy} P.
Arnoux, G.Rauzy ``Representation geometrique des suites the
complexite $2n + 1$''. В этой работе построена динамическая система
для слов с функцией роста $T(n)=2n+1$, обладающих дополнительным
комбинаторным свойством. В работе \label{Rote} \cite{R} G.Rote,
``Sequences with subword complexity 2n'' в терминах эволюции графов
Рози описаны слова с функцией роста $2n$.

Одним из основных результатов данной работы является обобщение этих
результатов на слова c линейной функцией роста, то есть с функцией
роста $T(n)=kn+l$, для $n>N$.

Перекладывания отрезков естественным образом служат обобщением
вращения круга. Эти преобразования были введены Оселедецом
 \cite{Os}, следовавшим идее Арнольда \cite{Arn}, (см. также \cite{KS}).
Рози \cite{Ra3} впервые показал, что связь между вращениями круга и
последовательностями Штурма обобщается если рассматривать
перекладывания отрезков. В связи с этим (в той же работе) он задал
вопрос описания последовательностей, связанных с перекладываниями
отрезков.

Такие последовательности являются еще одним естественным обобщением
слов Штурма. В частном случае, для $k=3$ отрезков, описание таких
последовательностей было получено в работе \cite{FHZ}, а работе
\cite{FZ2} были изучены частные случаи последовательностей,
порождаемых перекладыванием $4$-х отрезков. Стоит также отметить
работы \cite{AHP,Ba1,Ba2,BMP}.

В случае произвольного числа отрезков также получен ряд интересных
результатов. В работе \cite{FZ} получен комбинаторный критерий на
порождаемость слов, получаемых симметричным перекладыванием
отрезков, то есть перекладыванием, связанным с перестановкой $(1\to
k,2\to k-1,\ldots,k\to 1)$.

В работе \cite{FZ3}независимо от нас получен критерий порождаемости
слов преобразованием перекладывания отрезков, удовлетворяющих
следующему условию: траектория каждой концевой точки отрезка
перекладывания не попадает на концевую отрезка перекладывания, в том
числе сама на себя. В этом случае, как не сложно видеть, слова будут
иметь функцию сложности $T(n)=(k-1)n+1$. В данной работе получен
более общий результат, не требующий выполнения этого условия.

Отметим, что во всех этих работах изучаются перекладывания, не
меняющие ориентацию отрезков, а характеристические множества
совпадают с отрезками перекладывания. В нашей работе с помощью языка
графов Рози мы сначала изучаем слова, порождаемые
кусочно-непрерывным преобразованием отрезка, а затем показываем
эквивалентность множества таких слов и слов, порождаемых
перекладываниями. Этот метод дает возможность описать р.р. слова,
связанные с произвольным перекладыванием отрезка, более того, мы не
требуем, чтобы характеристические множества, соответсвующие символам
алфавита, совпадали с перекладываемыми отрезками.

\paragraph {Цель работы.} Данная диссертация посвящена
 исследованию нормальных базисов алгебр медленного роста,
а также исследованию взаимосвязи между комбинаторными свойствами
слов, цепочками порождающих их автоматов (графов Рози) и
порождающими их динамическими системами. Свойства периодичности
такаже находятся в центре внимания настоящей работы.

\paragraph{Научная новизна.}
Все основные результаты являются новыми. Среди них можно отметить
следующие:

\begin{enumerate}

\item Описание нормальных базисов
граничных алгебр, то есть алгебр  с функцией сложности,
асимптотически равной $n+C$.

\item Построение общего критерия порождаемости слова
преобразованием перекладывания отрезков в автоматных терминах,
решение вопроса, поставленного Рози.

\item Описание множества  сбалансированных
непериодических слов над произвольным алфавитом в терминах
порождающей динамической системы.

\end{enumerate}

\paragraph{Основные методы исследования.}

Основными инструментами исследования в работе являются исследование
цепочки автоматов (графов Рози), порождающих сверхслово, а также
техника подстановок, восходящих к А.~И.~Ширшову. Мы пользуемся
различными результатами теории равномерно-рекуррентных слов и слов
Штурма. Также используются результаты эргодической теории, такие как
существование инвариантной меры и равномерность иррациональных
сдвигов тора.

\paragraph{Практическая и теоретическая ценность.}

Работа носит теоретический характер. Результаты диссертации могут
быть полезны в комбинаторной теории колец и полугрупп, в частности,
при изучении мономиальных алгебр, а также в символической динамике.

\paragraph{Апробация результатов.}

Основные результаты диссертации докладывались на следующих
семинарах:
\begin{enumerate}
\item ``Кольца и модули'' кафедры высшей алгебры МГУ в 2000-2006
гг.

\item ``Арифметика и геометрия'' кафедры теории чисел МГУ в 2004-2005г.

\item Dynamics seminar, Einstein Institute of Mathematics под руководством профессора Х.Фюрстенберга (H.Furstenberg) в 2004г.

\item Семинаре под руководством профессора А.Френкеля (A.~Freankel) в 2004г., Weizman Institute of Science.

\item "Динамические системы" кафедры дифференциальных уравнений МГУ под руководством Ю.С. Ильяшенко в 2008 г.

\end{enumerate}

\paragraph{Публикации.}

Результаты диссертации опубликованы в следующих работах:
\begin{enumerate}

\item Белов А. Я. Чернятьев А.Л., Слова медленного роста и
перекладывания отрезков // Успехи Мат. Наук, 2008, 63:1(379),
159–160. В этой работе Белову А.Я. принадлежит постановка задачи и
формулировка основной теоремы, Чернятьеву А.Л. принадлежат
доказательства основных теорем 1 и 2.

\item Чернятьев А. Л.,Сбалансированные слова и символическая
динамика, Фундаментальная и прикладная математика, Том 13, выпуск 5,
2007 г., стр. 213-224.

\item Чернятьев А.Л.,  Белов А.Я., Описание слов Штурма над алфавитом
из n символов,  Математические методы и приложения, 6-й  мат. Симп.
МГСУ, - Москва, МГСУ, 1999г., стр. 122-128. А.Я. Белову принадлежит
постановка задачи и формулировка основной теоремы, Чернятьеву А.Л.
принадлежит построение конструкции динамической системы и
доказательство основной теоремы.

\item Белов А.Я, Чернятьев А.Л., Описание множества слов, порождаемых
перекладыванием отрезков // Депонировано в ВИНИТИ №1048-B2007. В
этой работе Белову А.Я. принадлежит постановка задачи, понятие
размеченного графа Рози и формулировка основной теоремы, Чернятьеву
А.Л. принадлежат доказательства основных теорем 1 и 2.

\item Чернятьев А.Л. Слова с минимальной функцией роста // Вестник
МГУ, 2008 г.

\end{enumerate}

\paragraph{Структура и объем работы.} Диссертация состоит из оглавления, введения, четырех глав и списка литературы, который включает 73 наименования.

\paragraph{Благодарности.}

Автор глубоко благодарен своем научным руководителям --- доктору
физико-математических наук Алексею Яковлевичу Белову и доктору
физико-математических наук, профессору Александру Васильевичу
Михалеву за постановку задач, обсуждение результатов и постоянное
внимание к работе.

Также автор хотел бы поблагодарить за внимание и обсуждения работы
доктора физико-математических наук, профессора Виктора Николаевича
Латышева, доктора физико-математических наук Николая Германовича
Мощевитина, доктора физико-математических наук Андрея Михайловича
Райгородского, доктора физико-математических наук, профессора Юлия Сергеевича Ильяшенко.

Автор выражает свою отдельную благодарность профессору университета
Вайсмана А.Френкелю (A.~Freankel) за детальное обсуждение работы, а
также профессору Л.Замбони ( Luca Q. Zamboni) за обсуждение $4$-ой главы работы.

\paragraph{Краткое содержание работы.}

{\bf Первая глава} посвящена обзору и доказательству базовых
результатов комбинаторики слов.

Отдельно рассматриваются результаты теории слов Штурма, представлено
доказательство теоремы эквивалентности, которое используется в
следующих главах. Также доказываются основные результаты из теории
графов Рози.

{\bf Вторая глава} посвящена изучению сбалансированных слов над
произвольным алфавитом. Изучаются свойства слов, порождаемых
иррациональными сдвигами тора.

Пусть $W$ -- слово над бинарным алфавитом, порождаемое сдвигом
одномерного тора, то есть динамикой $(\mathbb{S}^1,U,T_\alpha,x_0)$.

Обозначим $\beta =\alpha/q+r/q$, $U_q=\{x|qx\in U\}$, $y_0=x_0/q$,
где $q$, $r$~--- целые. Несложно показать, что динамика
$(\mathbb{S}^1,T_{\beta},U_q,y_0)$ порождает то же слово $W$.
Переход от первой динамики ко второй мы назовем {\it
$q$--размножением}

Доказана следующая теорема об эквивалентности динамик:

\begin{theorem}
Пусть два слова $W_1$ и $W_2$ порожденные динамиками
$(\mathbb{S}^1,T_{\alpha},U,x_0)$ и $(\mathbb{S}^1,T_{\beta},V,y_0)$
совпадают. Тогда существуют $p$ и $q$ такие, что множества $U$ и $V$
при соответствующих размножениях совпадают c точностью до
поворота,т.е. $U_p=T_{\delta}(V_q)$ для некоторого $\delta>0$.
\end{theorem}

Далее проводим редукцию от сбалансированного слова $W$ над
$n$-буквенным алфавитом к $n$ бинарным словам Штурма:

Построим слова $W_1,W_2, \ldots ,W_n$ над бинарными алфавитами

$$A_1=\{a_1,\bar{a_1}\}, A_2=\{a_2,\bar{a_2}\}, \ldots
,A_n=\{a_n,\bar{a_n}\}$$ следующим образом:

$W_i=(w^i_n)_{n\in\mathbb{Z}}$

$$
w^i_n=w_n=\left\{
\begin{array}{rcl}
   a_i,\  w_n=a_i\\
   \bar{a_i},\  w_n \neq a_i\\
\end{array}\right.
$$

Каждому слову $W_i$ соответствует одномерная динамика
$(\mathbb{S}^1,T_{\alpha_i},\Delta_i, x_i)$, ее порождающая,
следовательно слово $W$ порождается свигом на $n$-мерном торе с
вектором сдвига $\gamma=(\alpha_1,\ldots, \alpha_n)$.

Обозначим через $M$ замыкание траектории начальной точки
$x=~(x_1,x_2,\ldots, x_n)$ при действии $T_\gamma$.

Доказано следующее предложение:

\begin{proposition}
$M$ гомеоморфно множеству $\mathbb{S}^1\times \{1,2,\ldots ,N\}$.
\end{proposition}

Теперь видно, что динамика реализуется на множестве \par
$M=\mathbb{S}^1\times\{1,2,\ldots, N\}$, а отображение имеет
вид:\par $f:(x,k)\to (x+\alpha, k+1 \mod N)$. \par Будем теперь
понимать под $U_i$ характеристическое множество для символа $a_i$,
которое лежит на замыкании траектории, и пусть также $U^k_i$
обозначает часть характеристического множества, которое лежит на
$k$-ой компоненте связности (окружности) $M$. Далее изучаются
соответствующие характеристические множества.

Основным результатом данной главы является
\begin{theorem}
Пусть $W$ -- сбалансированное непериодическое слово над алфавитом
$A$. Тогда для $W$ существует динамическая система $(M,f)$,
удовлетворяющая следующим условиям:
\begin{enumerate}
\item $M=\mathbb{S}^1 \times {\mathbb Z}_m$ как топологическое
пространство.

\item $f:M \to M$ есть композиция поворота на $\alpha$ в
$\mathbb{S}^1$ и сдвига на $1$  в ${\mathbb Z}_m$. Длина
$\mathbb{S}^1$ равна $m$.

\item Каждая компонента $\mathbb{S} \times \{k\}$ $k=1, \ldots ,n$
разбита на $2m$ дуг: $m$ красных и $m$ синих; все красные имеют
длину $\alpha$, все синие $1-\alpha$,красные и синие дуги
чередуются.

\item Синий цвет имеет $l$ оттенков, красный -- $k$ оттенков,
$k+l=|A|$-- число букв в алфавите. Все середины дуг даного оттенка
образуют вершины правильного многоугольника (``правильный
$1$-угольник'' -- это точка на окружности,``правильный
$2$-угольник'' -- пара диаметрально противоположных точек).

\item При переходе от компоненты $\mathbb{S} \times \{k\}$ к
компоненте $\mathbb{S} \times \{k+1\}$ ($\mathbb{S} \times
\{m+1\}=\mathbb{S} \times \{1\}$) порядок расположения оттенков
внутри красных и синих компонент сохраняется, а сами красные и синие
дуги (``рулетки'') проворачиваются относительно друг друга на $1$.
Так что преобразование $f$ приводит к смещению на $\alpha$
относительно красных компонент и на $1-\alpha$ (в обратную сторону)
относительно синих.

\end{enumerate}
\end{theorem}

{\bf Третья глава} посвящена описанию нормальных базисов {\it
граничных} алгебр, алгебр, которые связаны с другим обобщенем слов
Штурма -- словам {\it медленного роста}.

Слово $W$ называется {\it словом медленного роста}, если
$F_W(n)=n+K$, для всех $n\geq N$. В словах Штурма $F_W(n)=n+1$ и
$F_W(n+1)-F(n)=1$. Естественным обобщением будут слова, для которых
$F(n)=n+K$, то есть выполняется $F(n+1)-F(n)=1$ для достаточно
больших $n$.

Показывается, что для каждого такого слова существует граф Рози,
соответствующий некоторому слову Штурма. Действительно, если начиная
с некоторого $n>N$ выполняется $T(n)=n+K$, то это означает, что для
каждого $n>N$ в слове $W$ существует только одно правое и одно левое
специальное слово длины $n$. То есть в графах Рози таких слов есть
одна входящая и одна выходящая развилка, а значит, существует слово
Штурма с такой же эволюцией $k$-графов, что и данное. Дальше
просходит редукция к теореме эквивалентности. Доказана следующая
\begin{theorem}\label{IntroMinGrow}
Пусть $W$ -- рекуррентное слово над произвольным конечным алфавитом
$A$. Тогда следующие условия на слово $W$ эквивалентны:
\begin{enumerate}

\item Существует такое натуральное $N$, что функция роста слова
$W$ равна $T_W(n)=n+K$, для $n\geq N$  и некоторого постоянного
натурального $K$.

\item Существуют такое иррациональное $\alpha$ и целые $n_1,n_2,
\ldots n_m$, что слово $W$ порождается динамической системой
$(\mathbb{S}^1,T_\alpha, I_{a_1},\ldots ,I_{a_n},x)$, где $T\alpha$
-- сдвиг окружности на иррациональную величину $\alpha$, $I_{a_i}$
-- объединение дуг вида $(n_j\alpha,n_{j+1}\alpha)$.
\end{enumerate}
\end{theorem}

Для произвольной алгебры $A$ через $V_A(n)$ обозначается размерность
вектрного пространства, порожденного мономами длины не больше $n$,
$V_A(n)$ называется {\it функцией роста} алгебры $A$. Пусть
$T_A(n)=V_A(n)-V_A(n-1)$. Если алгебра однородна, то $T_A(n)$ есть
размерность векторного пространства, порожденного мономами длины
ровно $n$. $T_A(n)$ называется {\it функцией сложности} алгберы $A$.

Известно (см. \cite{BBL}) что либо $\lim_{n\to
\infty}{(T_A(n)-n)}=-\infty$ (в этом случае есть альтернатива): либо
$\lim{V_A(n)}=C<\infty$ и тогда $\dim A < \infty$, либо
$V_A(n)=O(n)$ и алгебра имеет {\it медленный рост}), либо  $T_A(n)-n
< \Const$, либо, наконец, $\lim_{n\to \infty}{(T_A(n)-n)}=\infty$.

В последнем случае рост может быть хаотичным, поэтому для изучения
интересны первые два случая. Случай, когда $T_A(n)-n < \Const$ (т.е.
случай алгебр {\it медленного роста}) исследовался Дж.~Бергманом и
Л.~Смоллом. Нормальные базисы для таких алгебр исследованы в работе
\cite{BBL}. Назовем алгебру {\it граничной}, если $T_A(n)-n<\Const$.

Описание нормальных базисов граничных алгебр следует из теоремы
\ref{IntroMinGrow} и результатов \cite{BBL}.

{\bf Четвертая глава} посвящена описанию слов, связанных с
перекладыванием отрезков. Рассматриваются случаи, когда ориентация
отрезков сохраняется и когда нет. В терминах графов Рози дается
комбинаторное описание сверхслов порождающихся перекладыванием
отрезков. В центре внимания слова, имеющие функцию роста
$T_W(n)=kn+l$. Если слово обладает такой функцией роста, то
$T_W(n+1)-T_W(n)=k$.

Известно, что если перекладывание $k$ отрезков {\it регулярно}, то
есть траектория любого из концов отрезков перекладывания не попадает
на другой конец любого отрезка, то эволюция любой точки является
словом с ростом $T_W(n)=kn+l$.

Мы ищем условия на слово $W$, при которых оно порождалось бы
преобразованием перекладываний отрезка, не обязательно являющегося
регулярным. Отметим, что сдвиг окружности, по сути, является
перекладыванием двух отрезков с сохранением ориентации.

Рассмотрим соответствие между подсловами и подмножествами $M$. Легко
видеть, что если начальная точка принадлежит множеству $U_i$, то ее
эволюция начинается с символа $a_i$. Рассмотрим образы множеств
$U_i$ при отображениях $f^{(-1)}, f^{(-2)}, \ldots$. Ясно, что если
точка принадлежит множеству
$$T^{(-n)}(U_{i_n}) \cap T^{-(n-1)}(U_{i_{n-1}}) \cap \ldots \cap T^{(-1)}(U_{i_1}) \cap U_{i_0},$$
то эволюция начинается со слова $a_{i_0}a_{i_1}\cdots a_{i_n}$.
Соответственно, количество различных существенных эволюций длины
$n+1$ равно количеству разбиений множества $M$ на непустые
подмножества границами подмножств $\partial U_i$ и их образами при
отображениях $f^{-1}, f^{-2}, \ldots , f^{-n+1}$. Обозначим через
$I_u$ множество разбиения, которое соответствует слову $u$. Ясно,
что специальным подсловам соответствуют те интервалы, которые
делятся образами концов перекладываемых интервалов. Для данного
слова $u$ назовем слово $v$ {\it левым (соотв. правым) потомком},
если $u$ -- суффикс (соотв. префикс) слова $v$, в соответствии с
этим будем называть вершину в $G_n$ левым (соотв. правым) потомком
вершины в $G_k$, $n>k$. Прообраз конца интервала может являться
граничной точкой только для двух интервалов, соответственно,
специальные подслова могут иметь валентность только равную $2$.
Сформулируем

{\bf Правило 1.} Для того, чтобы бесконечное слово $W$ порождалось
системой $(I,T,U_1,\ldots,U_k)$ необходимо, чтобы любое специальное
слово имело валентность $2$.

Таким образом, мы можем наложить условие на эволюцию графов Рози:
начиная с некоторого $k$ все $k$-графы Рози имеют входящие и
исходящие развилки степени $2$. Предположим, что некоторому подслову
$w$ соответствует характеристический интервал, полностью лежащий
внутри интервала перекладывания. Пусть точка $A\in [0,1]$ делит
$I_w$ на два интервала, образы которых лежат в  $I_{a_k}$  и
$I_{a_l}$ соответственно, а точка $B\in [0,1]$ -- делит на
интервалы, прообразы которых лежат в $I_{a_i}$ и $I_{a_j}$
соответственно.

Выбор минимального невстречающегося  слова, а, значит, удаляемого
ребра, определяется взаиморасположением точек $A$ и $B$, а также
сохранением или сменой ориентации отображения на этих множествах.
Итого, имеется 8 вариантов, которые разбиваются на четыре пары,
соответствующие одинаковым наборам слов. Например, слову $a_i w a_k$
соответствует ситуация

$$
B<A, T^{-1}([x_w,B])\subset I_{a_i} , T([x_w,A]\subset I_{a_k} ).
$$

Введем понятие {\it размеченного графа Рози}. Граф Рози называется
{\it размеченным}, если

\begin{enumerate}
\item Ребра каждой развилки помечены символами $l$ (``left'') и $r$ (``right'')

\item Некоторые вершины помечены символом ``--''.
\end{enumerate}

{\it Последователем} размеченного графа Рози назовем ориентированный
граф, являющийся его последователем как графа Рози, разметка ребер
которого определяется по правилу:
\begin{enumerate}

\item Ребра, входящие в развилку должны быть помечены теми же символами, как и ребра, входящие в любого левого потомка этой вершины;

\item Ребра, выходящие из развилки должны быть помечены теми же символами, как и ребра, выходящие из любого правого потомка этой вершины;

\item Если вершина помечена знаком ``--'', то все ее правые потомки также должны быть помечены знаком ``--''.

\end{enumerate}

{\bf Замечание.} Поясним смысл разметки графа. Пусть ребра входящей
развилки соответствуют $a_i$ и $a_j$, символы $l$ и $r$
соответствуют левому и правому множеству в паре
$(T(I_{a_i}),T(I_{a_j}))$. Если символы $a_k$  и $a_l$ соответствуют
ребрам исходящей развилки, то символы $l$ и $r$ ставятся в
соответствии с порядком ``лево-право'' в паре $(I_{a_k},I_{a_l})$.
Знак ``--'' ставится в вершине, если характеристическое множество,
ей соответствующее, принадлежит интервалу перекладывания, на котором
меняется ориентация.

Условие для перехода от графа $G_n$ к $G_{n+1}$:

{\bf Правило 2.}
\begin{enumerate}
\item Если в графе нет двойных развилок, соответствующих биспециальным подсловам, то при переходе от $G_n$ к $G_{n+1}$ имеем $G_{n+1}=D(G_n)$;

\item Если вершина, соответствующая биспециальному слову не помечена знаком ``--'', то ребра, соответствующие запрещенным словам выбираются из пар $lr$ и $rl$
\item Если вершина помечена  знаком ``--'', то удаляемые ребра должны выбираться из пары $ll$ или $rr$.
\end{enumerate}

Назовем эволюцию размеченных графов Рози {\it правильной}, если {\bf
правила 1} и {\bf 2} выполняются для всей цепочки эволюции графов,
начиная с $G_1$, назовем эволюцию {\it асимптотически правильной},
если {\bf правила 1} и {\bf 2} выполняются, начиная с некоторого
$G_n$. Будем говорить, что эволюция размеченных графов Рози {\it
ориентированна}, если в $k$-графах нет вершин, помеченных знаком
``--''.

\begin{theorem}
Равномерно-рекуррентное слово $W$
 \begin{enumerate}
\item  Порождается перекладыванием отрезков, тогда и только тогда, когда слово обеспечивается асимптотически правильной эволюцией размеченных графов Рози.
\item  Порождается перекладыванием отрезков  с сохранением ориентации тогда и только тогда, когда слово обеспечивается асимптотически правильной ориентированной эволюцией размеченных графов Рози.
\end{enumerate}
\end{theorem}

%% file: chapter1.tex
\section{Пространство слов и символическая динамика.}\label{Chapter1}

\subsection{Пространство слов.}
В этой части мы определим основные понятия комбинаторики слов. В
дальнейшем $А$ будет обозначать конечный алфавит, то есть непустое
множество элементов (символов).  Через $A^+$ обозначим множество
всех конечных последовательностей, символов, или {\it слов}.

Конечное слово всегда может быть единственным образом
представлено в виде

$w = w_1 \cdots w_n$, где $w_i \in A, 1\leq i \leq n$. Число $n$
называется {\it длиной} слова $w$ и обозначается $|w|$

Множество $A^+$ всех конечных слов над $A$ образует простую
полугруппу, где полугрупповая операция определяется как
конкатенация (приписывание).

Если к множеству слов добавить элемент $\Lambda$ (пустое слово),
то получим свободный моноид $A^*$ над $A$. Длина $\Lambda$ по
определению равна 0.

Слово $u$ есть {\it подслово} (или {\it фактор}) слова $w$, если
существуют слова $p,q\in A^+$ такие, что $w = puq$.

Если слово $p$ (или $q$) равно $\Lambda$, то $u$ называется {\it
префиксом} (или {\it суффиксом}) слова $w$.

Мы будем обозначать $\Pref(w)$ ( $\Suf(w)$) множество всех
префиксов (соответственно суффиксов) слова $w$.

Для любой пары целых чисел таких, что $1 \leq i \leq j\leq n$ мы
обозначим через  $w[i,j]$ подслово $w[i,j] = w_i\cdots w_j$.

Если $w=w_1\cdots w_n$ $w_i\in A$, $i=l, \cdots n$ -- слово, то
{\it обратное слово} $\tilde{w}$ для $w$ определяется как
$\tilde{w} = w_n \cdots w_l$. Более того, мы полагаем, что
$\tilde{\Lambda}=\Lambda$.

Слово называется {\it палиндромом}, если $\tilde{w} = w$.
Множество всех палиндромов обозначается $\PAL$.

Пусть $w=w_1\cdots w_n$ $w_i\in A$, $i=l, \ldots ,n$ -- некоторое
слово. Натуральное число $q$ называется {\it периодом} $w$, если
$w_i=w_i+q$ для всех $i\in [1,\ldots ,n-q]$. Мы будем обозначать
$p_w$, или просто $p$, минимальный период $w$.

Слово $w$ называется {\it периодическим}, если $p\leq [|w|/2]$
(здесь $[x]$ обозначает целую часть $x$)  Через $\mathbb{N}$
(соответственно $\mathbb{N}_+$) мы будем обозначать неотрицательные
(положительные) целые числа.

{\it Одностороннее} (соответственно {\it двустороннее}) {\it
бесконечное слово} над алфавитом $A$ --  это отображение

$w:\mathbb{N}_+\to A$, (соответственно $w:\mathbb{Z} \to A$).

Для каждого $n$, мы полагаем $w_n=w(n)$ и обозначим
$W=w_1w_2\cdots$.

Слово $u\in A^+$ -- конечное подслово $W$, если существуют такие
$i,j\in N$, ($1\leq i\leq j$), что $u=w_i\cdots w_j$.
Последовательность $w[i,j]=w_i \cdots w_j$ назовем вхождением $u$
в $W$.

Слово $W$ назовем в {\it существенно периодическим}, если его
можно представить как $W=uv^{\infty}=uvvv\cdots$, где $u\in A^*$,
$v \in A^+$.

Бесконечное в обе стороны слово назовем просто {\it периодичным},
если оно имеет вид
$$
W= v^\infty= \cdots vvv \cdots
$$

Через $F(W)$ обозначим множество всех подслов (конечных и
бесконечных) слова $W$.

\begin{definition}
Два бесконечных слова $W$ и $V$ над алфавитом $A$ назовем {\bf
эквивалентными}, если $F(W)=F(V)$.
\end{definition}

{\it Морфизмами} на пространстве  конечных слов мы будем называть
обычные полугрупповые морфизмы, то есть отображение
$\varphi:A^*\to A^*$ -- морфизм слов, если оно сохраняет операцию
конкатенации (приписывания): $\varphi(uv)=\varphi(u)\varphi(v)$
для любых $u,v \in A^*$. Это определение естественным образом
продолжается для морфизмов бесконечных слов.

Понятия {\it мономорфизма,  эпиморфизма и изоморфизма} вводятся
обычным путем.

\subsection{Рекуррентность и равномерная рекуррентность}

\begin{definition}
Слово $W$ называется {\bf рекуррентным}, если каждое его подслово
встречается в нем бесконечно много раз (в случае двустороннего
бесконечного слова, каждое подслово встречается бесконечно много
раз в обоих направлениях). Слово $W$ называется {\it
равномерно-рекуррентным} или ({\it р.р словом}), если оно
рекуррентно и для каждого подслова $v$ существует натуральное
$N(v)$, такое, что для любого подслова $W$ $u$ длины не менее, чем
$N(v)$, $v$ является подсловом $u$.
\end{definition}

Пусть $W$ -- бесконечное слово. Для любого подслова $v$ можно
определить множество {\it возвращаемых} слов $\Ret_W(v)$, а
именно, слово $u$ -- {\it возвращаемое} для $v$, если $vuv$ --
подслово $W$ и $v$ -- не является подсловом $u$. Ясно, что для
рекуррентных слов множество возвращаемых слов $\Ret(v)$ каждого
подслова $v$ будет непустым, а в случае равномерной рекуррентности
множество длин слов из $\Ret(v)$ будет ограниченным.

\begin{definition}
Для произвольного слова $W$ можно определить функцию роста:
$$
T_W(n)=\Card F_n(W)
$$
\end{definition}
Ясно, что если $T_W(n)=0$ для какого-то $n$, то $W$ -- конечное
слово. В противном случае бесконечное.

\begin{theorem}                           \label{Th1:4}
Следующие два свойства бесконечного слова $W$ равносильны:

а) для любого $k$ можно найти $N(k)$ такое, что любой участок $W$
длины $k$ можно найти в любом участке $W$ длины $N(k)$;

б) если все конечные куски слова $V$ являются конечными кусками
слова $W$, то и все конечные куски слова $W$ являются конечными
кусками слова $V$. \Noproof
\end{theorem}

\begin{lemma}[о компактности]\index{Лемма о компактности}
                                        \label{Lemma11}%
Пусть ${\cal M}$~--- множество слов неограниченной длины над
конечным алфавитом ${\cal A}$. Тогда существует бесконечное слово
$V$, все подслова которого являются подсловами слов из ${\cal M}$.
\Noproof
\end{lemma}

\begin{corollary}
Если $u$~--- подслово слова $W$, и в слове $W$ есть сколь угодно
длинные подслова, не содержащие $u$, то существует слово $W'$,
которое беднее кусками, чем $W$. При этом можно выбрать $W$ так,
чтобы $u$ не входило в $W'$. \Noproof
\end{corollary}

\begin{corollary}
Для любой убывающей в смысле отношения $\sqsupseteq$ цепочки
сверхслов найдется инфимум. \Noproof
\end{corollary}

\begin{corollary}
Для слова $W$ найдется беднейшее по подсловам слово $W'$ такое,
что $W \sqsupseteq W'$. \Noproof
\end{corollary}

Мы имеем следующую теорему, означающую, что из кусков любого слова
можно составить равномерно-рекуррентное слово.

\begin{theorem}              \label{Th1:5} \label{ThUnifRec}
Пусть $W$~--- бесконечное слово. Тогда существует равномерно
рекуррентное слово $\widehat{W}$, все подслова которого являются
подсловами $W$. \Noproof
\end{theorem}

Эта теорема исключительно важна в комбинаторике слов, т.к. очень
часто позволяет сводить изучение произвольных слов к изучению р.р.
слов.

\subsection{Специальные подслова}

\begin{definition}
Рассмотрим бесконечное (одностороннее или двустороннее) слово над
алфавитом $A$. Пусть $v$ -- его подслово и $x\in A$. Тогда
\begin{enumerate}
\item Символ $x$ -- левое (правое) расширение $v$, если $xv$
(соотв. $vx$) принадлежит $F(W)$.

\item Подслово $v$ называется  левым (правым) специальным
подсловом, если для него существуют два или более левых (правых)
расширения.

\item Подслово $v$ называется биспециальным, если оно является и
левым, и правым специальнм подсловом одновременно.

\item Количество различных левых (правых) расширений подслова
назовем  левой (правой) валентностью этого подслова.

\end{enumerate}
\end{definition}

\subsection{Морфизмы}
Классическим методом получения бесконечных (в одну сторону) слов
является итерация {\it подстановок} (или {\it морфизмов}).
Рассмотрим отображение множества конечных подслов над алфавитом
$A$ в себя: $\varphi: A^*\to A^*$, которое сначала определим на
буквах алфавита:
$$
\begin{array}{ll}
a_1\to \varphi(a_1)\\
a_2\to \varphi(a_2)\\
\vdots \\
a_n\to \varphi(a_n).\\
\end{array}
$$

Продолжим это отображение на конечные слова:
$$
\varphi(x_1x_2\cdots x_k)=\varphi(x_1)\varphi(x_2)\cdots
\varphi(x_k)
$$.

Для бесконечных (вправо) слов отображение продолжается
естественным образом:
$$
\varphi(x_1 x_2
x_3\cdots)=\varphi(x_1)\varphi(x_2)\varphi(x_3)\cdots\\
$$
В случае, если для какого-то символа $a\in A$ $\varphi(a)$
начинается с $a$, то имеем последовательность префиксов:
$$
\begin{array}{ll}
   \varphi(a)=au, \ u\in A^*\\
   \varphi^2(a)=\varphi(au)=\varphi(a)\varphi(u)\\
   \varphi^3(a)=\varphi(\varphi^2(a))=\varphi^2(a)\varphi^2(u)\\
   \vdots\\
\end{array}
$$

И таким образом, последовательность
$a,\varphi(a),\varphi^2(a),\ldots$ задает бесконечное слово $W$
такое, что $\varphi(W)=W$. Такие слова называются {\it
инвариантными относительно подстановки $\varphi$ }.

{\bf Пример 1.} Слово Фиббоначи.\label{FibbSeq} Рассмотрим
последовательность слов над алвафитом $A=\{a,b\}$:\par

$f_0=b$, $f_1=a$,$f_2=ab$, $f_3=aba,\ldots$, $f_n=f_{n-1}f_{n-2}$.
Предел этой последовательности слов, т.е. слово
$W=abaababaabaab\cdots$, называется {\it словом Фиббоначи}. Оно
инвариантно относительно подстановки $\varphi(a)=ab$,
$\varphi(b)=a$

{\bf Пример 2.} Слово Туе-Морса. \label{Morse} Рассмотрим
подстановку $\varphi(a)=ab$, $\varphi(b)=ba$. Тогда:
$$
\begin{array}{ll}
   \varphi^2(a)=\varphi(ab)=abba\\
   \varphi^3(a)=\varphi(abba)=abbabaab\\
   \vdots\\
\end{array}
$$

Слово, инвариантное относительно $\varphi$, называется {\it словом
Туе-Морса (Thue-Morse word)}. Это слово обладает следующим
свойством: оно свободно от кубов, то есть не содержит подслов вида
$u^3$. Легко проверить, что не существует бесконечных слов над
двубуквенным, свободных от квадратов, то есть подслов вида $u^2$.

{\bf Пример 3.} Слово Триббоначи.\label{Tribb} Рассмотрим слово
над трехбуквенным алфавитом $A=\{a,b,c\}$, инвариантное
относительно подстановки $\varphi(a)=ab$, $\varphi(b)=ac$,
$\varphi(c)=a$, называется {\it словом Триббоначи}.

\subsection{Вопросы периодичности}
Самым простым примером бесконечного слова являются периодические,
или существенно периодические (с предпериодом) слова. Имеет место
следующая
\begin{theorem}
Пусть $W$ -- бесконечное вправо слово над конечным алфавитом $A$
$(\Card(A)\geq 2)$. Тогда следующие условия эквивалентны:
\begin{enumerate}
\item Слово $W$ существенно периодично.

\item В слове $W$ конечное число правых специальных подслов.

\item Существует $n$ такое, что для каждого подслова длины большей
$n$ существует ровно одно возвращаемое слово

\item Существует $n$ такое, что $T_W(n)=T_W(n+1)$
\end{enumerate}
\end{theorem}

\begin{proposition} Множество $L$ конечных слов является множеством подслов бесконечного
слова $W=w_1w_2w_3\cdots$ тогда и только тогда, когда

\begin{enumerate}
\item Каждое слово длины $n$ в $L$ содержится в некотором слове
длины $n+1$ в $L$ и

\item Подмножества $L_1,L_2,L_3, \ldots$, где $L_n ={w \in L :
|w|=n}$,
\end{enumerate}

\end{proposition}

\subsection{Графы Рози подслов.}

Удобным инструментом для описания слова $W$ является {\it графы
подслов}, или графы Рози (Rauzy's graphs), введенные Рози
\cite{Ra}, которые строятся следующим образом: {\it $k$-граф}
слова $W$ -- ориентированный граф, вершины которого
взаимнооднозначно соответствуют подсловам длины $k$ слова $W$, из
вершины $A$ в вершину $B$ ведет стрелка, если в $W$ есть подслово
длины $k+1$, у которого первые $k$ символов -- подслово
соответствующее $A$, последние $k$ символов -- подслово,
соответствующее $B$. таким образом, ребра $k$-графа биективно
соответствуют ($k+1$)-подсловам слова $W$.

\vbox{
\begin{center}
\includegraphics{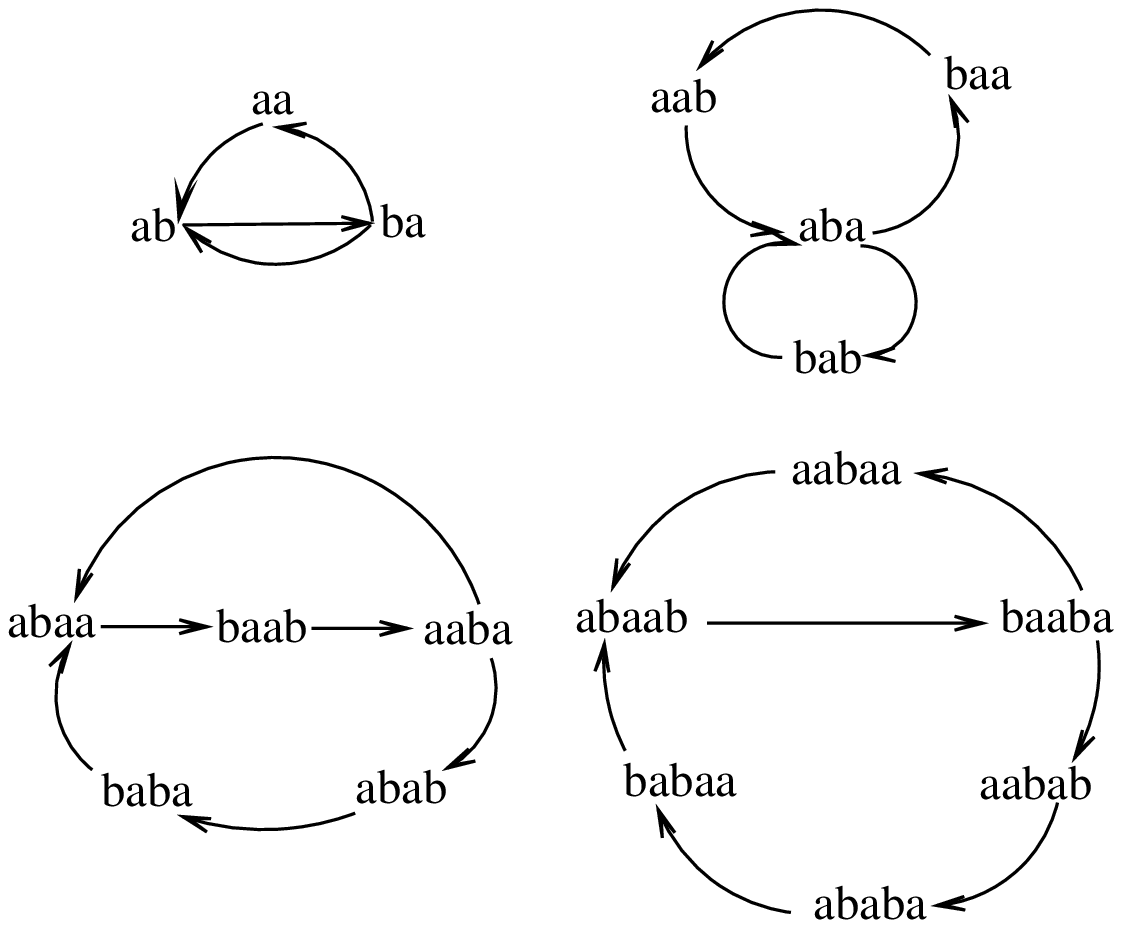}\\
{\it 2-,3-,4-,5-графы Рози для слова Фиббоначи.}
\end{center}
}\label{RisGraphRauzy1}

Ясно, в $k$-графе $G$ слова $W$ правым специальным словам
соответствуют вершины, из которых выходит (соотв. в которые
входит) больше одной стрелки. Такие вершины мы будем называть
развилками. Граф $G$ будем называть {\it сильно связным}, если из
любой вершины в любую вершину можно перейти по стрелкам.

{\it Последователем} ориентированного графа $G$ будем называть
ориентированный граф $\Fol(G)$ построенный следующим образом:
вершины графа $G$ биективно соответствуют ребрам графа $G$, из
вершины $A$ в вершину $B$ ведет стрелка, если в графе $G$ конечная
вершина ребра $A$ является начальной вершиной ребра $B$.

Связность графов Рози и рекуррентность соответсвующего слова
связаны естественным образом. Имеет место следующее
\begin{proposition}\label{reccur}
Пусть $W$ -- бесконечное (в одну сторону) слово. Следующие условия
эквивалентны:

\begin{enumerate}
 \item Слово $W$ -- рекуррентно.

 \item Для всех $k$ соответствующий
$k$-граф слова $W$ является сильно связным.

\item Каждое подслово $W$ встречается не меньше двух раз.

\item Любое подслово является продолжаемым слева.
\end{enumerate}
\end{proposition}

\Proof Достаточно доказать следующие 4 импликации:

$(i) \to (ii)$. Пусть $w$ и $v$ -- два подслова длины $k$ слова
$W$. Так как $W$ -- рекуррентно, то в $W$ существует подслово вида
$wAvBw$, где $A,B$ -- некоторые подслова. Соответственно, в $G_k$
существует путь из вершины, соответствующей $w$ в вершину $v$ и
путь из вершины, соответствующей $v$ в вершину, соответствующую $w$.
Значит, $G_k$ сильно связанный.

$(ii) \to (iii)$. Предположим противное. Пусть $v$ -- подслово,
встречающееся ровно раз в $W$, тогда существует префикс $v'$,
который содержит $v$ и который тоже встречается ровно раз в $W$.
Но тогда $G_k$ (где $k$ -- длина $v'$) не является строго
связанным.

$(iii) \to (iv)$. Пусть $v$ -- подслово $W$ и слово $v$
встречается в $W$ минимум дважды. Ясно, что $v$ является
продолжаемым (от второго вхождения к первому).

$(iv) \to (i)$. Предположим противное, пусть $W$ не рекуррентно.
Пусть $v$ подслово, которое встречается конечное число раз в $W$.
Тогда пусть $v'$ -- префикс $W$ достаточно большой длины, то есть
содержащий все вхождения $v$. Слово $v'$ больше не встречается в
$W$ следовательно, не является продолжаемым слева.

\Endproof

В терминах графов Рози можно выразить такие важные понятия, как
рост подслов, множество запрещенных подслов, минимальные
запрещенные слова и т.д.

\subsection{Слова, порождаемые динамическими системами.}

Пусть $M$ -- компактное метрическое пространство, $U\subset M$~--- его открытое подмножество, $f:M\to M$ -- гомеоморфизм компакта в
себя и $x\in M$ -- начальная точка.

По последовательности итераций можно построить бесконечное слово
над бинарным алфавитом:

$$
w_n=\left\{
\begin{array}{rcl}
   a,\ f^{(n)}(x_0)\in U\\
   b,\ f^{(n)}(x_0)\not\in U\\
\end{array}\right.
$$
которое называется {\it эволюцией} точки $x_0$. Символическая
динамика исследует взаимосвязь свойств динамической системы
$(M,f)$ и комбинаторных свойств слова $W_n$.

Для слов над алфавитом, состоящим из большего числа символов нужно
рассмотреть несколько характеристических множеств: $U_1,\ldots
,U_n$.

Заметим, что эволюция точки корректно определена только в случае,
когда траектория точки не попадает на границу характеристических
множеств $\partial U_1,\partial U_2,\ldots$.

Для того, чтобы рассматривать траекторию произвольной точки, мы
введем понятие {\it существенной эволюции}.

\begin{definition}
Конечное слово $v^f$ называется {\it существенной конечной
эволюцией} точки $x^{\ast}$, если в любой окрестности точки
$x^{\ast}$ существует открытое множество $V$, любая точка $x\in V$
из которого обладает эволюцией $v^f$.  Бесконечное слово $W$
называется {\it существенной эволюцией} точки $x^{\ast}$, если
любое его начальное подслово -- существенная конечная эволюция
точки $x^{\ast}$.
\end{definition}

Под {\it эволюцией точки}, когда это не вызывет недоразумений,
будем понимать, существенную эволюцию. Отметим, что точка может
иметь несколько существенных эволюций.

\begin{proposition}Пусть $V$ -- конечное слово. Тогда множество точек с конечной
существенной эволюцией $V$ замкнуто. Аналогичное утверждение верно
для бесконечного слова $W$.
\end{proposition}

\Proof См. \cite{BK}

\begin{definition}\label{MorphDin}
Морфизмом двух динамик $G:(M_1,f_1)\to (M_2,f_2)$ называется
непрерывное отображение $g:M_1\to M_2$ такое, что диаграмма

$$
\begin{array}{ccc}
M_1&\stackrel{f_1}{\longrightarrow} & M_2\cr
\downarrow\lefteqn{g}&&\downarrow\lefteqn{g}\cr
M_1&\stackrel{f_2}{\longrightarrow} & M_2\cr
\end{array}
$$
коммутативна.
\end{definition}

Понятие мономорфизма, эпиморфизма и изоморфизма вводится обычным
образом.

Фактор динамика естественным образом определяется на
фактор-топологии тогда и только тогда, когда $f$ переставляет
классы эквивалентности отображения $f$. Отметим, что прообразы
точек при морфизмах замкнуты.

\subsection{Функция рассогласования.}
Вместе с метрикой Хэмминга, разделяющей любые два различные слова
можно рассмотреть другую меру, называемую {\it плотностью
рассогласования}.
\begin{definition}

Пусть $V=\{v_n\}$ и $W\{w_n\}$ -- два бесконечных слова над
алфавитом $A$. Функция $\rho=\rho(V,W)$ называется плотностью
рассогласования и определяется по формуле:
$$
\rho(U,W)=\lim_{n\to \infty}{\frac{\sum_{i=-n}^n \rho(i)}{2n}}
$$

где
$$
\rho(i)=\left\{
\begin{array}{rcl}
   0,& v_i=w_i\\
   1,& v_i\neq w_i\\
\end{array}\right.
$$ есть функция рассогласования.

\end{definition}

\subsection{Соответствие между словами и разбиениями множества.}\label{Sootv:2}
Теперь посмотрим соответствие между словами и подмножествами $M$. Из
построения следует, что если начальная точка принадлежит множеству
$U_i$, то ее эволюция начинается с символа $a_i$. Рассмотрим образы
множеств $U_i$ при отображениях $f^{(-1)}, f^{(-2)}, \ldots$, $n\in
\mathbb{N}$.

Ясно, что если точка принадлежит множеству
$$f^{(-n)}(U_{i_n}) \cap f^{-(n-1)}(U_{i_{n-1}}) \cap \ldots \cap f^{(-1)}(U_{i_1}) \cap U_{i_0},$$
то эволюция начинается со слова $a_{i_0}a_{i_1}\cdots a_{i_n}$.

Соответственно, количество различных существенных эволюций длины
$n+1$ равно количеству разбиений множества $M$ на непустые
подмножества границами подмножств $\partial U_i$ и их образами при
отображениях $f^{(-1)}, f^{(-2)}, \ldots ,f^{(-n)}$.

{\bf Замечание.} Количество конечных существенных эволюций
напрямую связано с топологической размерностью множества $M$.
Ясно, например, что если $М$ гомеоморфно отрезку или окружности,
то одна точка может являться граничной только для двух открытых
подмножеств $M$ и, соответственно, иметь только две существенные
эволюции. Если $M$ гомеоморфно части плоскости $\mathbb{R}^2$ --
то существенных эволюций может быть сколь угодно много.

%\subsection{Частные случаи символических динамик.}\label{GenSym}

\subsection{Перекладывания отрезков.}\label{PerOtr}

Перекладывание отрезков является естественным обобщением сдвига
окружности -- в случае разбиения окружности на дуги длины $\alpha$
и $1-\alpha$ и величины сдвига, равной $\alpha$, это преобразование
совпадает с перекладыванием двух отрезков.

Более того, перекладывание отрезков является очень важным
преобразованием в эргодической теории, теории динамических систем
и теории чисел.

Рассмотрим общий случай:

Пусть отрезок $[0,1]$ разбит на полуинтервалы длин
$\lambda_1,\lambda_2,\ldots ,\lambda_n$ и $\sigma\in S_n$ --
перестановка на множестве $\{1,2,\ldots,n\}$.

Интервалы разбиения могут быть представлены через длины отрезков
следующим образом:

$$
X_i=\left [ \sum_{i<j} \lambda_j,\sum_{i\leq j} \lambda_j\right)
$$

Перекладывание отрезков ``переставляет'' отрезки $(X_1,X_2,\ldots
,X_n)$ разбиения между собой, в результате получается новое
разбиение

$$
(X_\sigma(1),X_\sigma(2), \ldots ,X_\sigma(n))
$$.

Более точно, преобразование $T$ ставит в соответствие каждой точке
$x\in X_i$ :
$$
T(x)=x+a_i
$$
где
$$
a_i=\sum_{k<\sigma^{-1}(i)}\lambda_{\sigma(k)}-\sum_{k<i}\lambda_k
$$

Если же преобразование переворачивает некоторый отрезок, то все
точки дополнительно симметрично отражаются относительно середины
этого отрезка. В случае, если перекладывание не переворачивает
отрезки, мы будем говорить, что такое перекладывание сохраняет
ориентацию.

\begin{definition}
Перекладывание отрезков $T$ называется {\it регулярным}, если для
любой точки $a_i$, где $X_i=[a_i,a_{i+1})$, $T^n(a_i)\neq a_j$

\end{definition}

Сформулируем следующую важную теорему:
\begin{theorem}
Перекладывание отрезков является регулярным тогда и только тогда,
когда траектория любой точки всюду плотна в $[0,1]$
\end{theorem}

\begin{theorem}
Пусть слово $W$ является эволюцией некоторой точки при регулярном
перекладывании $k$ отрезков. Тогда любое подслово $W$ имеет ровно
$k$ возвращаемых слов.
\end{theorem}

Изучение свойств слов, порождаемых перекладываниями отрезков
осуществляется теми же методами, что и при сдвиге
единичной окружности. Основным инструментом здесь является
рассмотрение отрицательных орбит концов интервалов перекладывания:

$$
0 = a_1 < a_2 < \ldots < a_{k+1} = 1,
$$

где

$X_i = [a_i, a_{i+1})$, ($i \in \{1, \ldots , k\}$).

Множество концов интервалов перекладывания $\{a_i| 1 \leq i \leq k
+ 1\}$ обозначим $X^{1}$. Слово  $w$ длины $n$ является подсловом
эволюции точки $x$, то есть бесконечного слова $U(x)$, тогда
только тогда, когда существует такой интервал $I_w \subset[0, 1]$
и точка $y \in I_w$ такие, что слово $w$ равно конкатенации
следующих символов:

$$
\mathcal{I}(x)\mathcal{I}(T(x))\cdots \mathcal{I}(T^{n-1}(x)) = w,
$$
где $\mathcal{I}(x)=a_i \in A $ тогда и только тогда, когда $x\in
X_i$

\begin{proposition} Пусть $T$ -- регулярное перекладывание $k$
отрезков. Тогда эволюция $U(x)$ произвольной точки $x$ имеет функцию сложности
$T_{U(x)}(n)=n(k-1)+1$, для любого $n\in\mathbb{N}$.
\end{proposition}

\subsection{Слова Штурма.} В этой части мы дадим  основные
определения и конструкции, связанные со словами Штурма (также
называемыми последовательностями Бетти, сбалансированными словами
и т.д.), а также сформулируем и докажем теорему эквивалентности.

\begin{definition}
Бесконечное вправо (соответственно, двусторонее бесконечное) слово
$W=(w)_{n\in \mathbb{N}}$ (соответственно, $W=(w)_{n\in
\mathbb{Z}}$ ) над конечным алфавитом $A$ называется словом
Штурма, если его функция сложности равна $T_W(n)=n+1$ для всех
$n\geq 0$.
\end{definition}

{\bf Примечание 1.} Так как $T_W(1)=2$, то слово Штурма по
определению является словом над двухсимвольным алфавитом. Здесь и
далее в этой главе мы будем рассматривать слова только над
бинарными алфавитами.

\begin{proposition}
Если $W=(w)_{n\in \mathbb{N}}$ -- слово Штурма, то $W$
рекуррентно.
\end{proposition}

\Proof Предположим противное и некоторое подслово $u$ встречается
в $W$ конечное число раз. Пусть оно не встречается, начиная с
$N$-ой позиции. Рассмотрим слово $V=(v)_{n\in \mathbb{N}}$:
$v_k=w_k+N$. Тогда все подслова $V$ являются подсловами $W$ и в
$V$ нет подслова $u$, значит $T_V(n)\leq n$ для какого-то $n$, а
следовательно, с какого-то момента периодично. \Endproof

{\bf Примечание 2.} В случае бесконечного в обе стороны слова при
использовании аналогичного определения слова Штурма возникают
исключительные классы слов, например:

$V_1=\cdots aaabbb\cdots$, $V_2= \cdots aaabaaa\cdots $. Поэтому в
определении двусторонних слов Штурма мы будем также предполагать
рекуррентность.

{\bf Пример.} Слово Фиббоначи.\label{FibbSeq} Рассмотрим
последовательность слов над алфавитом $A=\{a,b\}$:\par

$f_0=b$, $f_1=a$,$f_2=ab$, $f_3=aba,\ldots$, $f_n=f_{n-1}f_{n-2}$.
Напомним (см. \ref{FibbSeq}), что предел этой последовательности
слов, т.е. слово $W=abaababaabaab\cdots$, называется {\it словом
Фиббоначи}.

Для подслов небольшой длины нетрудно заметить:
$$
\begin{array}{ll}
T_W(1)=\Card\{a,b\}=2\\
T_W(2)=\Card\{aa,ab,ba\}=3,\\
T_W(3)=\Card\{aba,baa,bab,aab\}=4,\\
T_W(4)=\Card\{abaa,abab,baab,baba,aaba\}=5\\
\end{array}
$$

\begin{proposition}
Слово Фиббоначи является словом Штурма.
\end{proposition}
\Noproof

\subsection{Структура слов Штурма.}
\begin{proposition}
Слово Штурма содержит ровно одно из двух подслов: $aa$ или $bb$.
\end{proposition}
\Proof Действительно, $T_W(2)=3$, однако, если слово содержит и
$aa$ и $bb$, то оно обязано, в силу рекуррентности, содержать и
$ab$ и $ba$, что приводит к противоречию. \Noproof

\begin{definition}
Мы будем говорить, что слово Штурма имеет тип $"a"$, если в нем
встречается подслово $aa$, и типа $"b"$, если встречается подслово
$bb$.
\end{definition}

Теперь перейдем к определению слов Штурма через другое понятие --
свойство сбалансированности.

\begin{definition}
Слово $W$ называется сбалансированным, если для любых его двух
подслов $u,v$ одинаковой длины выполняется неравенство:

\begin{equation} \label{bal}
||u|_a-|v|_a| \leq 1
\end{equation}

\end{definition}

Ясно, что из неравенства (\ref{bal}) автоматически следует
$||u|_b-|v|_b|\leq 1$

Естественным обобщением понятия сбалансированности является {\it
m-~сбалансированность}:

\begin{definition}
Слово $W$ над алфавитом $A$ называется $m$-сбалансированным по
символу $a\in A$, если для любых подслов $u,v\subset W$ таких что
$|u|=|v|$ выполняется $||u|_a-|v|_a|\leq m$. Слово называется
$m$-сбалансированным, если оно сбалансировано по каждому символу
алфавита.
\end{definition}

Свойство сбалансированности означает, что для непериодического
слова $W$ символы $a$ и $b$ перемешаны  ``наилучшим образом''.
Если для какого-то $k$ будет выполняться более сильное условие, то
есть для любых двух слов длины $k$ количество символов $a$ (а
соответственно и $b$) было бы одинаковым, то слово $W$ было бы
периодичным с периодом $k$.

Докажем теперь следующее важное

\begin{proposition}\label{Prop:1-1}
Пусть слово $W$ не является сбалансированным. Тогда существует
такое слово $u$ (возможно, пустое), что слова $aua$ и $bub$
являются подсловами $W$.
\end{proposition}
\Proof Если слово $W$ не сбалансировано, то существуют два
подслова одинаковой длины $v_1$ и $v_2$ такие, что
$||v_1|_a-|v_2|_a|\geq 2$. Будем также считать, что эти слова
имеют наименьшую длину среди всех подслов, обладающих данным
свойством.

Обозначим  $v^k_1, v^k_2$ -- префиксы длины $k$ слов $v_1$ и $v_2$
.

Пусть также $d_k=||v^k_1|_a-|v^k_2|_a|$. Ясно, что $d_n\geq 2$
(где $n=|v_1|=|v_2|$) и $|d_{k+1}-d_k|$ равно $0$ или $1$.

Положив для удобства $d_0=0$ и рассматривая последовательность
$d_k$ от $k=0$ до $k=n$, мы обнаружим первый момент, когда
$d_k=2$. В силу минимальности длины слов $v_1$ и $v_2$, имеем
$k=n$

Пусть теперь $v_1=a_1 a_2 \cdots a_n$ и $v_2=b_1b_2 \cdots b_k$.
Тогда ясно, что $a_n\neq b_n$, иначе при обрезании $n$-х символов
мы опять получили бы более короткие слова с этим же свойством.
Аналогично, $a_1\neq b_1$. Также ясно, что $a_1=a_n$, $b_1=b_n$,
так как в противном случае слова $a_2\cdots a_{n-1}$ и $b_2\cdots
b_{n-1}$ также обладали бы этим свойством. Пусть $a_1=a_n=a$ и
$b_1=b_n=b$.

Но тогда $a_2=b_2$, в противном случае опять получатся более
короткие слова с разницей символов, б\'ольшей $1$.

 Рассматривая также далее последовательность префиксов,
индуктивно приходим к равенствам $a_k=b_k$, $k=2,\ldots ,n-1$,
следовательно, слово $u=a_2\cdots a_{n-1}$ -- искомое. \Endproof

Следующее утверждение является первой частью общей теоремы
эквивалентности:

\begin{theorem} \label{SturmBalance}
Слово $W$ является словом Штурма тогда и только тогда, когда оно
 существенно непериодично и сбалансированно.

\end{theorem}

\Proof Докажем, сначала {\bf импликацию $\Rightarrow$}.
Предположим противное, то есть слово $W$ не является словом
Штурма. Покажем тогда, что оно не сбалансированно. Рассмотрим
такое наименьшее $n_0$, что $T_W(n_0)\geq n_0+2$. Ясно, что
$n_0\geq 2$, так как $T_W(1)=2$ по определению. Так как $n_0$ --
минимальное, то $T_W(n_0-1)=n_0+1$, это значит, что существует как
минимум два подслова $u$ и $v$, которые имеют по два разных
продолжения справа. Из минимальности $n_0$ опять же следует, что
$u$ и $v$ различаются только в первых символах, иначе более
короткие суффиксы также обладали бы этим свойством, следовательно
существует подслово $w$ такое, что $u=aw$ и $v=bw$. Так как $u,v$
имеют по два разных продолжения, то в $W$ встретятся подслова
$awa$ и $bwb$, значит $W$ не сбалансированно.

Теперь докажем {\bf импликацию $\Leftarrow$}. Предположим , что
$W$ не сбалансированно. Тогда по предложению \ref{Prop:1-1},
найдется такое подслово $u=w_0w_1 \cdots w_n$, что $aua$ и $bub$
-- также подслова $W$. Предположим также, что $u$ имеет
минимальную длину из всех слов с таким свойством.

Отметим, что, исходя из этого предположения, следует, что если
пара слов одинаковой длины $u_1,u_2$ не сбалансирована, то слова
имеют длину, не меньше, чем $n+3$. Действительно, из
доказательства предложения \ref{Prop:1-1}, следует, что
минимальная длина такого $u$ не больше наименьшей длины не
сбалансированной пары уменьшенной на 2. Поэтому любая
не сбалансированная пара имеет длину не меньше $|u|+2=n+3$.

Заметим, что по предложению $u$ не пусто, иначе в слове $W$
встретились бы слова $aa$ и $bb$.

По этой же причине, $w_0=w_n$, $w_1=w_{n-1}$ и более обще,
$w_k=w_{n-k}$, то есть $u$--{\it палиндром.}

Теперь, мы знаем, что различных слов длины $n+1$ ровно $n+2$.
Слово $u$ может быть расширено двумя путями как влево, так и
вправо. Соответственно, ровно одно из слов $au$ или $bu$ (мы
предположим, что $au$) может быть расширено вправо двумя путями.
Следовательно, $aua$, $aub$ и $bub$ -- подслова $W$, а $bua$ -- не
является подсловом. Пусть $i$ -- первый момент, когда подслово
$bub$ встречается в $W$. Докажем следующее
\begin{proposition}
Слово $au$ не является подсловом слова
$$
v=w_iw_{i+1}\cdots w_{i+2n+3}.
$$
\end{proposition}

\Proof Длина $v$ равна $2n+4$, длина $bub$ равна $n+3$ и длина
$au$ равна $n+2$, следовательно, в предложении в точности
утверждается, что первый символ вхождения $au$ не может
встречаться в $bub$. Предположим противное, то есть что начало
$au$ перекрывается с $bub$. Ясно, что перекрытие происходит не в
первой позиции $bub$, если перекрытие происходит позиции $w_k$, то
это означает, что $aw_0\cdots w_{n-k}=w_k \cdots w_nb$ и это
означает, что $w_k=a$ и $w_{n-k}=b$, что противоречит тому, что
$u$ -- палиндром.
\Endproof

{\bf Завершение доказательства Теоремы \ref{SturmBalance}}

Понятно, что в слове $v=w_i w_{i+1} \cdots w_{i+2n+3}$ встречается
ровно $n+3$ подслова длины $n+2$. Но из предложения следует, что
среди них нет подслова $au$, следовательно, какое-то подслово
встречается как минимум два раза. Все подслова из $v$ могут быть
продолжены вправо единственным образом, так как двумя способами
продолжается вправо только $au$.  Следовательно, сверхслово $W$
существенно периодично.
\Endproof

 Теперь отметим еще одну важную характеристику слов Штурма --
 плотность вхождения символа.

 Сформулируем следующее предложение -- определение.

\begin{proposition}
Определим {\it плотность вхождения символа $a$ } в слове
$W=(w)_{\mathbb{N}}$ как предел
$$
\rho(a)=\lim_{n\to \infty}|w_0 \cdots w_n|_a/n
$$
Если $W$ -- слово Штурма, то предел существует и $\rho(a)$
иррационально.
\end{proposition}

\Proof Сначала докажем, что этот предел существует, если $W$ --
слово Штурма. Пусть $a_n$ -- минимальное количество символов в
подлслове $W$ длины $n$. Из сбалансированности $W$ следует, что
максимальное количество вхождений равно $a_n+1$.  Ясно, что
достаточно доказать существование предела последовательности
$a_n/n$ и его иррациональность.

Слово длины $kq+r$ может быть разбито на $k$ слов длины $q$ и одно
слово длины $r$, из чего немедленно следует, что
\begin{equation} \label{Eq1:2}
ka_q \leq a_{kq+r}\leq k (a_q+1)+r
\end{equation}

Если рассмотреть достаточно большие $n$, то есть $n> q^2$, мы
можем записать $n=kq+r$, где $k \geq q$ и $0 \leq r < q$.  Так как
$r < k$, то нетрудно проверить, что $r/n < k/n < 1/q$,
следовательно, из второй части неравенства (\ref{Eq1:2}) имеем :
$a_n/n \leq a_q/q+2/q$.

Далее, так как $n \geq a_n \geq k a_q > r a_q$, то из первой части
неравенства (1) вытекает $(a_q-1)n/q=ka_q+(ra_q-n)/q \leq a_n$
значит, после деления на $n$ имеем:
$$
\frac{a_q}{q}-\frac{1}{q}\leq \frac{a_n}{n}\leq
\frac{a_q}{q}+\frac{2}{q}.
$$
Значит, последовательность $a_n/n$ является последовательностью
Коши, а, значит, имеет предел.

Предположим теперь, что этот предел является рациональным числом
$p/q$. Из неравенства $ka_n\leq a_kn < a_kn+1 \leq k(a_n+1)$,
следует, что если $n'$ делит $n$, то $a_n/n \leq a_{n'}/{n'}<
(a_{n'}+1)/n' \leq (a_n+1)/n$.

В частности, последовательность
$(a_{2^nq}/{2^nq})_{n\in\mathbb{N}}$ является возрастающей, а
$(a_{2^nq}+1/{2^nq})_{n\in\mathbb{N}}$ -- убывающей.

Тогда имеем неравенства:
$$
\frac{a_q}{q} \leq
\frac{a_{2^nq}}{2^nq}<\frac{a_{2^nq}+1}{2^nq}\leq \frac{a_q+1}{q}.
$$

Но эти последовательности стремятся к $p/q$, что возможно только
при $a_q=p$ и $a_{2^nq}=2^np$ для всех $n$

Другой замечательной характеристикой слов Штурма является их
описание через динамику на единичной окружности.

\begin{proposition} Пусть
$\mathbb{S}^1$ -- окружность единичной длины, $U\subset\mathbb{S}^1$
-- дуга длины $\alpha$, $T_\alpha$ -- сдвиг окружности на
иррациональную величину $\alpha$. Тогда динамика произвольной точки
порождает слово с функцией сложности $T_W(n)=n+1$, т.е. слово
Штурма.
\end{proposition}

{\bf Набросок доказательства (подробнее см. \cite{L1,L2}).}
Рассмотрим образы $U$ при итерациях $T^{(-1)}(U),T^{(-2)}(U),
\cdots$. Поскольку длина дуги совпадает с величиной сдвига, то при
первой такой итерации появится ровно одна новая граничная точка и,
соответственно, одно новое множество разбиения, при каждой следующей
итерации также будет появляться ровно одна дополнительная граничная
точка, остальные граничные точки будут переходить друг в друга. Так
как количество множеств разбиения -- это и есть количество подслов
данной длины, а изначально было два множества, то функция сложности
$T_W(n)=n+1$

Слова, получаемые такими динамическими системами, иногда называют
{\it механическими словами} ({\it mechanical words}).

Возникает естественный

{\bf Вопрос:} {\it Любое ли слово Штурма может быть получено из
сдвига окружности}.

Ответ на него положительный. Таким образом, можно сформулировать
следующую {\bf теорему эквивалентности}:

\begin{theorem}\label{TheorEq}
Следующие условия на слово $W$ эквивалентны:
\begin{enumerate}
\item Слово $W$ имеет функцию сложности $T_W(n)=n+1$.

\item Слово $W$ сбалансированно и не периодично.

\item Слово $W$ порождается системой $(\mathbb{S}^1, U, T_\alpha)$
с иррациональным $\alpha$.
\end{enumerate}
\end{theorem}
\Noproof

Рассмотрим подстановку над алфавитом $A=\{a,b\}$: $\sigma(a)=ab$ и
$\sigma(b)=a$. Тогда последовательность
$a,\sigma(a),\sigma^2(a).\ldots$ задает {\it слово Фиббоначи}(см.
\ref{FibbSeq}):
$$
W=abaababaabaab\cdots
$$

Подстановка $\sigma$ называется {\it морфизмом Фиббоначи}, так как
длина $\sigma^n(a)$ равна $n$-му числу Фиббоначи. Действительно,
легко проверить, что подстановка удовлетворяет рекуррентному
соотношению $\sigma^{n+2}(a)=\sigma^{n+1}(a)\sigma^{n}(a)$, а
значит, и рекуррентному соотношению на длины.

Подстановка $\sigma$ обладает еще одним замечательным свойством:
если слово $W$ над бинарным алфавитом является словом Штурма, то
$\sigma(W)$ также является словом Штурма.

Дадим теперь более общее
\begin{definition}
Морфизм $\varphi:A^*\to A^*$ называется {\bf морфизмом Штурма},
если для произвольного слова Штурма W слово $\varphi(W)$ также
является словом Штурма.
\end{definition}

Отметим, что свойство быть морфизмом Штурма является локальным, то
есть его можно проверить на конечных словах. Справедлива следующая

\begin{theorem}[S.Berstel, P.S\'e\'ebold]
Морфизм $\varphi$ является морфизмом Штурма тогда и только тогда,
когда образ слова $u=abbabbababbaba$ является сбалансированным.
\end{theorem}

Для примера, морфизм Фиббоначи переводит слово $u$ в слово \par
$\sigma(u)=abaaabaaabaabaaabaab$, которое является
сбалансированным.

Теперь введем следующее
\begin{definition}
Слово Штурма $W$ называется {\it характеристическим (или
стандартным)}, если все его префиксы являются левыми специальными
подсловами.
\end{definition}

Согласно теореме эквивалентности, для заданного $\alpha$
существует бесконечно много слов Штурма с плотностью символов
$\alpha$ и $1-\alpha$ соответственно, а именно слова, порожденные
динамиками $(\mathbb{S}^1,T_\alpha,\Delta,x_0)$, которые
отличаются выбором начальной точки $x_0$. Стандартное слово Штурма
отвечает начальной точке $x_0=0$.

Как показано в \cite{L1}, стандартные слова Штурма могут быть
получены как предел последовательности $s_n$ конечных слов,
задаваемых рекуррентным соотношением: $s_{n}=s^{d_n}_{n-1}
s_{n-2}$, где $n\geq 1$, $s_{-1}=b$, $s_0=a$ и последовательность
целых $d_n$ удовлетворяет условию $d_0\geq 0$, $d_n>0$, при $n>0$.
Последовательность $(d_n)$ связана с величиной сдвига $\alpha$
следующим образом: непрерывная дробь $[0,1+d_0,d_1,d_2,\ldots]$
равна в точности $\alpha$.

Еще одной характеристикой слов Штурма является их описание через
количество подслов-палиндромов. Напомним, что конечное слово
$w=w_1w_2\cdots w_n$ называется {\it палиндромом}, если
$w=\tilde{w}=w_nw_{n-1}\cdots w_1$.

{\it Палиндромическим замыканием} конечного слова $u$
 называется наименьшее по длине слово -- $u^{+}$, которое
 имеет $u$ в качестве префикса.

 {\bf Пример}. Для слова $u=abaabab$ палиндромическое замыкание
 равно $u^{+}=abaababaaba$.

Следующая характеристика стандартных слов Штурма получена A. de
Luca \cite{AdL}:
\begin{proposition}\label{StPal}
Бесконечное вправо слово $W$ над бинарным алфавитом $A=\{a,b\}$
является характеристическим  словом Штурма тогда и только тогда,
когда существует бесконечное слово $\Delta(W)=a_0a_1\cdots$, в
котором оба символа встречаются бесконечно много раз такое, что:

$$
W=\lim_{n\to\infty}u_n,
$$

где $u_0=\varepsilon$ и $u_{n+1}=(u_na_n)^{(+)}$, $n\geq 1$.
\end{proposition}

Более общая теорема была доказана A.Droubay  и J.Pirillo
(см.\cite{DP}):
\begin{proposition}
Бесконечное слово $W$ является словом Штурма тогда и только тогда,
когда для любого натурального $n$ слово $W$ имеет ровно одно
подслово-палиндром длины $n$, если $n$ -- четное и ровно два
подслова-палиндрома длины $n$, если $n$ -- четное.

\end{proposition}

Напомним, что для рекуррентного слова $W$ и его подслова $u$ {\it
возвращаемые слова} -- это максимальные по длине слова,
находящиеся между двумя последовательными вхождениями $u$ в $W$.

Для слов Штурма имеет место следующее характеристическое свойство:

\begin{proposition}
Равномерно рекуррентное слово $W$ является словом Штурма, если и
только если каждый его префикс имеет ровно два возвращаемых слова.
\end{proposition}

\subsection{Слова Арно-Рози.}
Естественными обобщениями слов Штурма являются сбалансированные
слова над произвольным алфавитом, и слова с медленным ростом, о
которых пойдет речь в главах $3$ и $4$. Другой подход к обобщению
слов Штурма был исследован в работах П. Арно и Г.Рози (см.
\cite{Ra2}), полученный класс слов носит название {\it слов
Арно-Рози (Arnoux-Rauzy words)} или {\it AR-слов}.

Бесконечное слово $W$ над $k$-буквенным алфавитом называется {\it
AR- словом}, если для каждого натурального $n\geq 1$ существует
ровно одно левое и ровно одно правое специальное подслово длины
$n$, причем валентность этих специальных подслов равна в точности
$k$.

Ясно, что слова Штурма являются AR-словами. Примером AR-слова над
алфавитом более чем из двух символов является слово Трибоначи
(\ref{Tribb}). Оно строится над алфавитом из $3$-х символов.

Графы Рози для $AR$-слов над $k$-буквенным алфавитом имеют вид.
схожий с графами Рози для слов Штурма. Они имеют ровно одну
вершину в которую входят $k$ ребер и ровно одну, из которой
выходят $k$ ребер. В случае, если соответствующее вершине
специальное подслово является палиндромом, то оно является
биспециальным, то есть и правым и левым одновременно. Остальные
вершины графа имеют одно входящее и одно исходящее ребро. Все
$AR$-слова также являются равномерно-рекуррентными словами.

$AR$-слово $W$ называется {\it характеристическим}, если все его
префиксы являются левыми специальными подсловами. В частности,
слово Трибоначчи является характеристическим словом.

Характеристические $AR$-слова могут быть описаны с помощью так
называемого {\it правила Рози (Rauzy's rules)}. Правило
заключается в способе получения бесконечного (право) слова
посредством итерации последовательности морфизмов определенного
вида. Точнее, для заданного алфавита $A$ рассмотрим семейство
морфизмов $\tau_a:A^*\to A^*$, $a\in A$:

$$
\tau_a(b)=\left\{
\begin{array}{rcl}
   a,& b=a\\
   ab,& b\neq a\\
\end{array}\right.
$$

Для заданного слова $w=a_1\cdots a_n$ определим
$\tau_w=\tau_{a_1}\cdots \tau_{a_n}$. Например,
$\tau_{abc}(a)=\tau_{ab}(ca)=\tau_{a}(bcba)=abacaba$.

Имеет место следующее
\begin{proposition}
Произвольное характеристическое AR-слово над алфавитом A может
быть получено как предел последовательности
($(\tau_{d_n}(a))_{n\geq 1}$ для некоторого $a\in A$,  где $d_n$
есть префикс длины $n$ некоторого бесконечного слова $\Delta$,
которое содержит бесконечно много вхождений каждого символа из
$A$.
\end{proposition}

Следующий результат является обобщением теоремы \ref{StPal} для
характеристических $AR$-слов был получен \cite{DJP}:

\begin{theorem}
Произвольное характеристическое $AR$-слово может быть получено как
предел последовательности $(u_n)_{n\geq 0}$ слов, где
$u_{n+1}=(u_na_n)^{+}$ -- палиндромическое замыкание $u_n a_n$,
$u_0=\varepsilon$, и бесконечное слово $\Delta=a_0a_1a_2\cdots$
содержит бесконечно много вхождений каждого символа из $A$.
\end{theorem}

%\subsection{Символические динамики, порожденные подстановками}.

%% file: chapter2.tex
\section{Сбалансированные слова и динамические системы}\label{balance}

\subsection{Введение и постановка проблемы.}
Понятие сбалансированности играет существенную роль при изучении
комбинаторных свойств слов, а также в теории биллиардов и теории
динамических систем.

Для слов над алфавитом $A=\{a_1,a_2, \ldots,a_n\}$ свойство
сбалансированности (или ``равномерной перемешанности'') означает,
что в любых двух подсловах одинаковой длины количество символов
одного сорта (например, символов $a_1$) почти одно и то же. Более
точно, для любого $1\leq i \leq n$, и любых двух подслов
одинаковой длины, количество символов $a_i$ в них отличается не
более, чем на $1$.

Бесконечные непериодические слова над алфавитом из двух символов
-- это в точности слова Штурма.

Из теоремы эквивалентности \ref{TheorEq} следует, что
сбалансированное непериодическое слово порождается сдвигом на
единичной окружности, то есть может быть получено следующим
образом:
$$
w_n=\left\{
\begin{array}{rcl}
   a,&T_\alpha^{(n)}(x_0)\in U\\
   b,&T_\alpha^{(n)}(x_0)\not\in U,\\
\end{array}\right.
$$

где $T_\alpha: x \to x+\alpha (\bmod 1)$ -- сдвиг на величину
$\alpha$, $U$ -- дуга длины $\alpha$ и $x_0$ -- начальная точка.

{\bf Наша цель} -- обобщение данного результата на сбалансированные
слова над более чем двухбуквенным алфавитом, то есть построение
для каждого сбалансированного слова такой динамической системы,
которая порождала бы данное слово.

Отметим, что описание сбалансированных слов в арифметической форме
были получены Р.Грэхемом \cite{GR} , а позднее описание в форме
бинарных слов Штурма было получено П.Хубертом \cite{H}.

Также отметим, что описание периодических слов Штурма связано с
известной проблемой Френкеля (Fraenkel conjecture, \cite{T2}) о
покрытии множества целых чисел последовательностями Бетти (другое
название последовательностей Штурма).

\subsection{Основные конструкции и определения.}

Пусть $A=\{a_1,a_2, \ldots,a_n\}$ -- конечный алфавит, $W$ --
бесконечное (в обе стороны) слово, то есть произвольная
последовательность символов из $A$.

Напомним, что {\it функция сложности} $T_W(n)$ определяется как
количество различных подслов длины $n$ слова $W$

{\it Слово Штурма} -- это бесконечное слово с функцией сложности
$T_W(n)$=n+1, $n \geq 1$.

Для конечного подслова $u\subset W$ обозначим величины $|u|$ --
длина слова и $|u|_{a_i}$ -- количество вхождений в подслово $u$
символов $a_i$. Ясно, что $|u|=\sum_{a_i\in A}|u|_{a_i}$.

\begin{definition}
Бесконечное слово $W$ над алфавитом $A=\{a_1,a_2, \ldots,a_n\}$
называется $m$-сбалансированным, если для любых двух подслов
$u,v\subset W$ и любого символа $a_i\in A$ выполняется неравенство
$$
||u|_{a_i}-|v|_{a_i}|\leq m.
$$
В случае, когда слово $1$-сбалансированно, будем просто говорить,
что оно сбалансированно.

\end{definition}

Отметим, что $0$-сбалансированность ведет к периодичности:

\begin{proposition}
Пусть для бесконечного слова $W$ существует такое натуральное $k$,
что для любых двух подслов $u,v$ длины $k$ и любого символа $a_i$
выполняется $|u|_{a_i}=|v|_{a_i}$. Тогда слово $W$ периодично
\end{proposition}
\Proof Фактически, сбалансированность означает минимальное
отклонение от периодичности, также, как и минимальный рост.

\subsection{Конечные сбалансированные слова}
\begin{proposition} Пусть $k\geq 2$, тогда для любого конечного слова $u$, $u^2$ сбалансированно тогда и
только тогда, когда сбалансированно $u^k$.
\end{proposition}

\Proof {\bf1.} Ясно, что если $u^2$ не сбалансированно, то и $u^k$
не сбалансированно.\par {\bf 2.} Пусть $u^k$  не сбалансированно,
тогда пусть $v$ и $v'$ -- два подслова $u^k$ минимальной длины
такие, что $||v|_a-|v'|_a| \geq 2$.

Если  $|v|\geq |u|$, то $v=xy$, $v'=x'y'$ и $|x|=|x'|=|u|$. Слова
$x$ и $x'$  сопряжены с $u$ и, значит, $|x|_a = |x'|_a = |u|_a$.
Из этого следует, что $||y|_a-|y'|_a| \geq 2$, что приводит к
противоречию с минимальностью $|v|$. Следовательно, $|v|< |u|$ и
слово $v$ является подсловом $u^2$. \Endproof

Для любого конечного сбалансированного слова $w$ над алфавитом
$A=\{a,b\}$, существует такое целое $n$, что между двумя
последовательными вхождениями символа $b$ находится ровно $n$ или
$n+1$ символ $a$.

Более того, каждое такое сбалансированное слово $w$ может
начинаться и оканчиваться не более, чем $n+1$ символом $a$ и,
значит, может быть представлено в виде $w = SA_1A_2\cdots A_kP$,
где $k\geq 0$, $S=a^l$ ($0\leq l \leq n+1$), $A_i \in \{ba^n,
ba^{n+1}\}$  ($1 \leq i \leq k$) и $P=ba^l$ ($0\leq l \leq n+1$).
Если слово $w$ (не обязательно сбалансированное) может быть
представлено в таком виде, мы будем говорить, что оно является
{\it $n$-сводимым} или, просто, {\it сводимым}.

Для заданного сводимого слова $w = SA_1A_2\cdots A_kP$ мы можем
построить новое слово, которое обозначим $\Red(w) = sa_1a_2\cdots
a_kp$ такое, что:

\begin{itemize}
\item   $s=\varepsilon$, если $S\neq a^{n+1}$,  иначе $s=b$

\item $a_i=a$, если $A_i=ba^n$

\item $a_i=b$, если$A_i=ba^{n+1}$

\item  $p=\varepsilon$, если $P\neq ba^{n+1}$, иначе $s=b$
\end{itemize}
Если $n\geq 1$, то операция сведения нетривиальна,
действительно, тогда $|\Red(w)|\leq |w|/2$ (так как  $|s| \leq |S|/
2$ , $p \leq |P|/ 2$ и $|a_i|  \leq|A_i|/ 2$ ).

С предложенной операцией сведения связана следующая

\begin{lemma}
Пусть конечное слово $w$ над алфавитом $A$  содержит не менее
одного символа $b$. Тогда следующие утверждения эквивалентны:
\begin{enumerate}

\item Слово $w$ сбалансированно.

\item Слово $w$  $n$-сводимо (для некоторого целого $n > 0$) и
слово $\Red(w)$ сбалансированно.
\end{enumerate}
\end{lemma}

\Proof Так как любое сбалансированное слово сводимо, достаточно
показать, что $n$-сводимое слово не сбалансированно, если  и
только если не сбалансированно слово $\Red(w)$. Сначала
предположим, что $w$ не сбалансированно. Из предложения
\ref{balance} следует, что существует такое подслово $v$, что
$ava$ и $bvb$ -- подслова $w$.

Если $|v|_b=0$, то $v=a^n$ или $v=a^{n+1}$. В этом случае
$a^{n+2}$ является подсловом $ava$ и тогда $w$ не является
$n$-сводимым.

Тогда существует такое целое $l$, что  $a^lb$ является прeфиксом
слова $v$. Так как $bvb$ -- подслово $w$, то $l = n$ или $l = n +
1$. Так как $a^{n+2}$ не является подсловом $w$ и $aa^l$ подслово
$av$, мы имеем $l = n$.

Аналогично, $ba^n$ является суффиксом $v$, значит $v =
a^nh^n(u)ba^n$, где $u\in \{a; b\}^*$. Так как $ba^nh^n(u)ba^nb$ --
подслово $w$, то $aua$ -- подслово $\Red(w)$.

Так как $a^{n+1}h^n(u)ba^{n+1}$ -- посдлово $w$, то $bub$ --
подслово  $\Red(w)$. Следовательно, $\Red(w)$ не сбалансированно.

Теперь, предположим, что $\Red(w)$ не сбалансированно.  Из
предложения \ref{balance}, существует такое подслово $v$, что
$ava$ и  $bvb$ -- подслова $\Red(w)$. Так как $bvb$ -- подслово
$\Red(w)$, то $a^{n+1}h^n(v)ba^{n+1}$ -- подслово $w$. Теперь
заметим, что если $ava$ -- суффикс $\Red(w)$, то существует целое
$l$, ($0 \leq l \leq n$) такое, что $h^n(ava)ba^l$ -- суффикс $w$.
Тогда $ba^nh^n(v)ba^nb$ является подсловом $w$, и, значит, $w$ не
сбалансированно.

\subsection{Расширение сбалансированных слов.}

В этой части мы покажем, что любое конечное сбалансированное слово
может быть продолжено до более длинного сбалансированного слова.
Последовательность сбалансированных слов, каждое из которых
является префиксом последующего, задает бесконечное (вправо)
сбалансированное слово. Имеет место следующее

\begin{proposition}

Для любого сбалансированного слова над алфавитом $A=\{a,b\}$
существует не менее двух символов $x$ и $y$ из $A$ таких, что $wx$
и $yw$ сбалансированы.

\end{proposition}

\Proof

Предположим, что $wa$ и $wb$ не сбалансированы, а $w$ --
сбалансированно. Тогда существуют такие подслова $v$ и $v'$, что

$ava$ -- суффикс $wa$, $bv'b$ -- суффикс $wb$, а $bvb$ и $av'a$ --
подслова $w$. Ясно также, что $|v|\neq |v'|$. Если $|v| < |v'|$,
то $v' = uav$ для некоторого слова  $u$. Тогда $ava$ и $bvb$
являются подсловами $w$ и $w$ не сбалансированно: противоречие.
Аналогично доказывается, что в случае $|v|>|v'|$, $av'a$ и $bv'b$
-- подслова $w$.

Тепреь, если $w$ сбалансированно, то и обратное слово
$\widetilde{w}$ тоже сбалансированно. Для него существует $y\in
A$, что $\widetilde{w}y$ -- сбалансированно, значит, и $yw$ --
сбалансированно.

\Endproof

Для непериодических сбалансированных слов над бинарным алфавитом
существует наглядная геометрическая интерпретация:  координатную
плоскость и произвольную прямую с иррациональным тангенсом
наклона. Рассмотрим такие точки на этой прямой, у которых хотя бы
одна координата -- целое число (отметим, что в силу
иррациональности тангенса наклона,точка с двумя целыми
координатами может быть только одна). Точкам с целой координатой
по оси $X$ сопоставим символ $a$, точкам с целой координатой по
$Y$ сопоставим $b$. Если точка имеет две целые координаты -- ей
сопоставляется $a$, если прямая имеет угол наклона меньше $\pi/4$,
и $b$ в противном случае. Слово, которое можно прочитать вдоль
прямой будет словом Штурма и сбалансированным непериодическим.

Эту конструкцию можно обобщить в $n$-мерный аналог. В этом случае
образуются классы т.н. кубических ($3$-мерный аналог) и
гиперкубические ($n$-мерный аналог) {\it биллиардные} слова.

Еще одной важной характеристикой бесконечных слов является {\it
плотность}:

\begin{definition}
Пусть $W=(w_n)_{n\in\mathbb{Z}}$ -- бесконечное слово над алфавитом
$A=\{a_1, a_2,\ldots,a_n\}$. Плотностью символа $a_i$ называется
предел (если он существует):
$$
\rho(a_i)=\lim_{n\to\infty}|w_{-n} \cdots w_{-1}w_0w_1 \cdots
w_n|_{a_i}/2n
$$

В случае бесконечного вправо слова:
$$
\rho(a_i)=\lim_{n\to\infty}|{w_0w_1\cdots w_n}|_{a_i}/n
$$

\end{definition}
 Если слово $W$~-- слово Штурма, то оно порождается сдвигом
окружности. Сдвиг на иррациональную величину является равномерным
(см. \cite{W}), поэтому плотность каждого символа совпадает с
мерой Хаара характеристического множества. То же верно и для
любого слова, порождаемого сдвигом окружности.

Случай периодических сбалансированных слов был изучен Тайдеманом
(\cite{T1}). Он показал, что любое периодическое сбалансированное
слово обладает плотностью по каждому символу.

\subsection{Слова, порождаемые повороотом окружности.}
В этой части мы сформулируем и докажем несколько полезных
утверждений о сдвиге окружности и словах, порождаемых сдвигом
окружности на иррациональную величину.

\begin{proposition}
\begin{enumerate}

{\bf 1.} Траектория произвольной точки при иррациональном $\alpha$
всюду плотна (лемма Кронекера)

{\bf 2.} Для любых двух точек слова, порождаемые этой динамикой
(т.е. существенные эволюции) эквивалентны (то есть множества их
конечных подслов совпадают).

\end{enumerate}
\end{proposition}

\begin{proposition}
Пусть $W=(w_n)$~--- слово, которое порождается динамикой
$(\mathbb{S}^1,T_{\alpha},U,x_0)$ и пусть $\beta
=\alpha/q+{r}/{q}$, $U_q=\{x|qx\in U\}$, $y_0={x_0}/{q}$, где $q$,
$r$~--- целые. Тогда эволюция точки $x_0$, порожденная динамикой
$(\mathbb{S}^1,T_{\beta},U_q,y_0)$ будет совпадать с $W$.
\end{proposition}

\Proof Непосредственно проверяется, что $T^{(n)}_\alpha(x_0)\in U$
тогда и только тогда, когда $T^{(n)}_\beta(y_0) \in U_q$. Такую
операцию над динамикой мы будем называть {\it $q$-размножением.}
Ясно также, что если $U$ -- симметричное подмножество, то есть при
сдвиге на величину $1/q$ ($q$-целое) переходит в себя, то аналогично
можно рассмотреть множество $U^q=\{qx|x\in U\}$ и сдвиг, в $q$ раз
больший. Такую операцию мы будем называть {\it $q$-сжатием}.

\begin{proposition}\label{qraz}
Пусть два слова $W_1$ и $W_2$ порожденные динамиками
$(\mathbb{S}^1,T_{\alpha},U,x_0)$ и
$(\mathbb{S}^1,T_{\beta},V,y_0)$ совпадают. Тогда существуют $p$ и
$q$ такие, что множества $U$ и $V$ при соответствующих
размножениях совпадают c точностью до поворота,т.е.
$U_p=T_{\delta}(V_q)$ для некоторого $\delta$.
\end{proposition}

\Proof Покажем сначала, что  $1$, $\alpha$ и $\beta$ линейно
зависимы над $\mathbb{Q}$ (т.е. существуют целые $m,n,k$ такие,
что $m\alpha+n\beta=k$). Действительно, предположим противное.
Пусть $\alpha$ и $\beta$ линейно независимы с единицей над
$\mathbb{Q}$ . Рассмотрим сдвиг на двумерном торе $\mathbb{T}^2$
$h:(x,y)\to(x+\alpha,y+\beta) \bmod 1$. Согласно лемме Вейля,
траектория любой точки при таком сдвиге всюду плотна. Тогда
существует $k$ такое, что $h^k(x_0,y_0)\in
U\times(\mathbb{S}^1/V)$, но это означает, что в $k$-ом символе
последовательности различаются. Противоречие. Значит
$\alpha={p/q}\beta+{r}/{s}$, где $p,q,r,s$--целые. Покажем, что
данные $p$ и $q$~--- искомые. Положим
${\alpha}_p={\alpha}/{p}={\beta}/q+{r}/{sp}$,
$\beta_q={\beta}/{q}$, тогда $ps{\alpha_p}=ps{\beta_q}=\gamma
 (\bmod 1)$. Совместим начальные точки сдвигом, пусть  величина
сдвига равна $\delta$. Рассмотрим сдвиг на $\gamma$. Поскольку
$\gamma$~--- иррационально, то траектория любой точки при сдвиге на
$\gamma$ всюду плотна. Поэтому, если $U_p\ne V_q$, то существует
момент, когда траектория точки попадет в $U_p\bigtriangleup V_q$
(здесь $A\bigtriangleup B$ обозначает симметрическую разность
множеств $A$ и $B$), т.е. слова в данной позиции различаются.

\Endproof

Это предложение можно усилить  и потребовать не полное совпадение
слов, а их совпадение на множестве позиций, имеющем положительную
меру в $\mathbb{Z}$, то есть потребовать ненулевую плотность
рассогласования.

Теперь, пусть окружность разбита на $k$ дуг:
$$
\{I_1=[\beta_1,\beta_2),I_2=[\beta_2,\beta_3), \ldots
,I_n=[\beta_n,\beta_{n+1})\},\ \beta_{n+1}=\beta_1.
$$

Для произвольного символа $a\in A$ мы также будем обозначать
характеристическое множество $I_a$

Рассмотрим сдвиг окружности на иррациональную величину $\alpha$ и
построим соответствующее слово $W$ над алфавитом
$A=\{a_1,a_2,\ldots,a_n\}$.

\begin{proposition}
Конечное слово $u=x_0 x_1 x_2 \cdots x_k$ является подсловом $W$
тогда и только тогда, когда
$$
\bigcap_{i=0}^k T_\alpha^{-i}(I_{x_i})\neq \emptyset
$$
\end{proposition}
\Proof Это утверждение является следствием леммы Кронекера.
Действительно, если пересечение не пусто, то из всюдуплотности
траектории следует, что траектория попадет в пересечение (конечное
пересечение полуинтервалов -- интервал или отрезок или
полуинтервал). Кусок эволюции, начиная с момента попадания
траектории в пересечение, с длиной,  равной $k+1$ будет совпадать
слову $u$.

\Endproof

\subsection{Сбалансированные слова над алфавитом из n символов.}

Пусть $W=(w_n)_{n\in \mathbb{Z}}$~ -- бесконечное в обе стороны
непериодическое сбалансированное слово над алфавитом $A=\{a_1,a_2,
\ldots ,a_n\}$.

Построим слова $W_1,W_2, \ldots ,W_n$ над бинарными алфавитами
\par $A_1=\{a_1,\bar{a_1}\},A_2=\{a_2,\bar{a_2}\}, \ldots
,A_n=\{a_n,\bar{a_n}\}$ \par следующим образом:

$$
W_i=(w^i_n)_{n\in\mathbb{Z}},
$$

$$
{w_n}^i=w_n=\left\{
\begin{array}{rcl}
   a_i,w_n=a_i\\
   \bar{a_i},w_n \neq a_i\\
\end{array}\right.
$$

\begin{proposition}
Слово $W_i$ являются сбалансированным над алфавитом $A_i$ для
любого $i$.

\end{proposition}
\Proof Из сбалансированности $W$ следует сбалансированность $W_i$
по символу $a_i$, действительно, если $u=w_jw_{j+1}\cdots
w_k\subset W$ и $u'=w^i_j w^i_{j+1} \cdots w^i_{k}$~ --
соответствующее подслово $W^i$, то $|u|_{a_i}=|u'|_{a_i}$.
Сбалансированность по одному символу для слов над бинарным
алфавитом влечет сбалансированность в целом. \Endproof

Согласно теореме эквивалентности \ref{TheorEq}, для каждого слова
$W_i$ существует динамика $(S^1,T_{\alpha_i},\Delta_i,x_i)$
($i=1,2,\ldots ,n$), которая порождает данное слово. Заметим, что
для любого $k$ среди символов $w^1_k, w^2_k,\ldots ,w^n_k$ должен
быть ровно один символ из алфавита $A=\{a_1,a_2,a_3,\ldots\}$ и
ровно $n-1$ символ из алфавита $\bar A\{\bar a_1,\bar a_2,\bar
a_3,\ldots\}$.

\begin{proposition} $1,\alpha_i,\alpha_j$
линейно зависимы над $\mathbb Q$ для любых $i$ и $j$.
\end{proposition}

\Proof Предположим противное и рассмотрим сдвиг на вектор
$(\alpha_i,\alpha_j)$ на двумерном торе
$\mathbb{T}^2=\mathbb{S}^1\times \mathbb{S}^1$. Согласно лемме
Кронекера-Вейля, траектория любой точки линейно независимы с $1$
над $\mathbb{Q}$. Тогда плотна и траектория точки $(x_i,x_j)$,
следовательно, на каком-то шаге она попадет в множество $U=\{(x,y|
x\in \Delta_i,y \in \Delta_j)\}$. Это означает, что на $n$-м месте
в слове $W$ одновременно стоят символы $a_i$ и $a_j$, что
невозможно. Противоречие. \Endproof

Фактичеcки, мы можем представить динамическую систему, которая
порождает слово $W$ следующим образом:
$$
(\mathbb{T}^n, T_\gamma, U_1, U_2,\ldots, U_n, (x_1,x_2,\ldots,
x_n)),
$$
где $T_\gamma:\mathbb{T}^n\to \mathbb{T}^n$ -- сдвиг на $n$-мерном
торе на вектор $\gamma=(\alpha_1,\alpha_2,\ldots ,\alpha_n)$, а
$U_i=\{(y_1,y_2,\ldots ,y_n)\in \mathbb{T}^n| y_i\in \Delta_i\}$.

Однако ясно, что в силу попарной рациональной зависимости над
$\mathbb{Q}$, динамика реализуется на меньшем множестве.

Обозначим через $M$ замыкание траектории точки $x=(x_1,x_2,\ldots,
x_n)$ при действии $T_\gamma$.

\begin{proposition}
$M$ гомеоморфно множеству $\mathbb{S}^1\times \{1,2,\ldots N\}$.
\end{proposition}

\Proof Так как все $\alpha_i$ попарно линейно зависимы с единицей
над $\mathbb{Q}$, то существует такое $\alpha$ (например,
$\alpha_1$ ) такое, что

$$
\begin{array}{lcl}
\alpha_1=r_1/s_1\alpha+p_1/q_1\\
\alpha_2=r_2/s_2\alpha+p_2/q_2\\
\vdots\\
\alpha_1=r_n/s_n\alpha+p_n/q_n,\\
\end{array}
$$

где $r_i,s_i,p_i,q_i\in \mathbb{Z}$.

То есть $\gamma=\alpha(r_1/s_1, \ldots,r_n/s_n)+(p_1/q_1,\ldots,
p_n/q_n)=\alpha\vec{\rho_1}+\vec{\rho_2}$, где
$\vec{\rho_1},\vec{\rho_2}$~-- рациональные векторы.

Отсюда видно, что, так как сдвиги между собой коммутируют, то
замыкание траектории есть прямая направленная вдоль вектора
$\vec{\rho_1}$ и все ее образы при сдвиге на вектора вида
$m\vec{\rho_2}$ ($m\in\mathbb{Z}$). Поскольку $\vec{\rho_1}$~ --
рациональный вектор, то прямая на торе будет замкнута, то есть
гомеоморфна окружности, а поскольку $\vec{\rho_2}$~-- рационален,
то различных образов этой окружности при сдвигах на
$m\vec{\rho_2}$ будет конечное число. \Endproof

Теперь видно, что динамика реализуется на множестве
$$
M=\mathbb{S} \times\{1,2,\ldots ,N\},
$$
а отображение имеет вид:
$$
f:(x,k)\to (x+\alpha, k+1 \bmod N).
$$
Будем теперь понимать под
$U_i$ характеристическое множество для символа $a_i$, которое лежит
на замыкании траектории, и пусть также $U^k_i$ обозначает часть
характеристического множества, которое лежит на $k$-ой компоненте
связности (окружности) $M$.

Из всюдуплотности траектории на $M$ следует, что каждая точка
принадлежит какому-то множеству $U^k_i$. Для составления общей
картины нам необходимо описание характеристических множеств.

\begin{proposition} \label{SetsOfCirc}
Каждое множество $U^l_i$ является объединением $k_i$ дуг,
периодически разбивающих окружность (то есть инвариантное
относительно сдвига на $1/k_i$), причем $k_i$ не зависит от номера
компоненты $l$.
\end{proposition}

\Proof Рассмотрим слово, порождаемое на $l$-ой компоненте и
характеристическом множестве $U^l_i$ которое мы построим следующим
образом: рассмотрим произвольную точку траектории $x_0$, которая
попадает на $l$-ю компоненту, через $N$ шагов преобразования $f$
она снова попадет на эту же окружность, сместившись на вектор
$N\alpha$, следующая итерация -- такое же смещение на $N\alpha$ и
т.д. При попадании в $U^l_i$ записываем $a_i$, при непопадании --
$\bar a_i$. Легко видеть, что это слово совпадает со некоторым
словом вида $\{v_k=w^i_{kN}\}$ составленного из $N$-х символов
слова $W_i$, то есть совпадает со словом, порожденным динамикой
$(\mathbb{S}^1,T_N\alpha_i,\Delta)$. Из предложения и следствия к
нему следует, что $U^k_i$ есть объединение $k_i$ дуг, периодически
разбивающих окружность.\Endproof

Будем обозначать длины дуг, которые образуют $U_i$ через
$\lambda_i$.

\begin{proposition}
Для каждого $i$ выполняется $\alpha=\lambda_i+l_i/k_i$ или
$\alpha=-\lambda_i+l_i/k_i$

\end{proposition}
\Proof Предложение также является следствием предложения
{\ref{qraz}}. \Endproof

 Теперь разобьем характеристические множества на два класса:
 те $U_i$, для которых $\alpha=\lambda_i+{l_i}/{k_i}$ для удобства
 окрасим в красный цвет (а различать различные множества будем по оттенкам),
 остальные аналогично окрасим в синий.
Заметим, что в силу данного соотношения, если из левого конца
красной сдвинуться на $\alpha$ то попадем в правый конец дуги того
же цвета (и оттенка), а если из правого конца дуги синего цвета
сместиться на $\alpha$, то попадем в левый конец дуги того же
цвета (левый и правый концы определяются согласно ориентации
окружности)

\begin{proposition}
Все дуги одного цвета имеют одинаковую длину.\\
\end{proposition}
\Proof Покажем для красного цвета. Так как
$\alpha=\lambda_i+l_i/k_i=\lambda_j+l_j/k_j$, то
$\lambda_i-\lambda_j\in \mathbb{Q}$, что возможно только когда
$\lambda_i=\lambda_j$. Для синего цвета действуем аналогично.
\Endproof

\begin{proposition}
На одной компоненте дуги, имеющие общую границу, имеют разный
цвет.
\end{proposition}

\Proof Предположим противное. Пусть какие-то две дуги $(x_1,x_2),
(x_2,x_3)$ -- красные. Предположим также, что другая соседняя с
дугой $(x_1,x_2)$, дуга $(x_0,x_1)$ -- синяя (случай двух синих дуг
рассматривается аналогично).  Рассмотрим образы точек $x_1,x_2$ при
сдвиге на $\alpha$. Так как $\alpha=\lambda+l_i/k_i$ для дуг
грасного цвета, то $x_1+\alpha$ является концом дуги того же цвета,
что и $(x_1,x_2)$, а $x_2+\alpha$ -- концом дуги того же цвета, что
и $(x_2,x_3)$. Так как длины красных дуг равны $\lambda$, то точка
$x_1+\alpha$ является левым концом дуги того же цвета что и
$(x_2,x_3)$.

Но $x_1$ также является правым концом дуги синего цвета,
следовательно, $x_1+\alpha$ должен быть левым концом дуги этого же
цвета, значит, $x_1$ одновременно является левым концом двух
разных дуг. Противоречие.\Endproof

Итак, мы получили, что на одной компоненте разбиение на
характеристические множества, или раскраска, выглядят следующим
образом:
\begin{enumerate}
\item Каждое характеристическое множество $U_i$ представляет собой
объединение дуг одной длины, на каждой компоненте середины дуг
образуют правильный $k_i$-угольник (то есть инвариантны при сдвиге
на $1/k_i$).

\item Дуги разбиты на два цвета, красный и синий. Дуги одного
цвета имеют одинаковую длину

\item На одной компоненте (окружности) красные и синие дуги
чередуются.
\end{enumerate}

Осталось понять как меняется раскраска при переходе от $i$-ой к
($i+1$)-ой компоненте. Пусть $m$ -- количество красных (а,
значит, и синих дуг на одной компоненте)

\begin{proposition}
При переходе от $i$-ой к $(i+1)$-ой компоненте дуги красного цвета
остаются на месте, дуги синего цвета поворачиваются на $1/m$.
Количество компонент $N$ кратно $m$.
\end{proposition}

\subsection{Основная теорема о непериодических сбалансированных словах над произвольным алфавитом} \label{OsnTh}

Для более наглядного представления основной теоремы будем считать,
что окружность на, которой реализуется динамика имеет длину $m$

\begin{theorem}
Пусть $W$-- сбалансированное непериодическое слово над алфавитом
$A$. Тогда для $W$ существует динамическая система $(M,f)$,
удовлетворяющая следующим условиям:
\begin{enumerate}
\item $M=\mathbb{S} \times {\Bbb Z}_m$ как топологическое
пространство.

\item $f:M \to M$ есть композиция поворота на $\alpha$ в
$\mathbb{S}$ и сдвига на $1$  в ${\mathbb{Z}}_m$. Длина
$\mathbb{S}^1$ равна $m$.

\item Каждая компонента $\mathbb{S}^1 \times \{k\}$ $k=1, \ldots ,n$
разбита на $2m$ дуг: $m$ красных и $m$ синих; все красные имеют
длину $\alpha$, все синие $1-\alpha$,красные и синие дуги
чередуются.

\item Синий цвет имеет $l$ оттенков, красный -- $k$, $k+l=|A|$ --
число букв в алфавите. Все середины дуг данного оттенка образуют
вершины правильного многоугольника (``правильный $1$-угольник'' --
это точка на окружности,``правильный $2$-угольник'' -- пара
диаметрально противоположных точки).

\item При переходе от компоненты $\mathbb{S}^1 \times \{k\}$ к
компоненте $\mathbb{S}^1 \times \{k+1\}$ ( $\mathbb{S}^1 \times
\{m+1\}=\mathbb{S}^1 \times \{1\}$) порядок расположения оттенков
внутри красных и синих компонент сохраняется, а сами красные и
синие дуги ("рулетки") проворачиваются относительно друг друга на
$1$. Так что преобразование $f$ приводит к смещению на $\alpha$
 относительно красных компонент и на $1-\alpha$ ( в обратную сторону)
 относительно синих.

\end{enumerate}
\end{theorem}

{\bf Пример.}

\vbox{
\begin{center}
\includegraphics{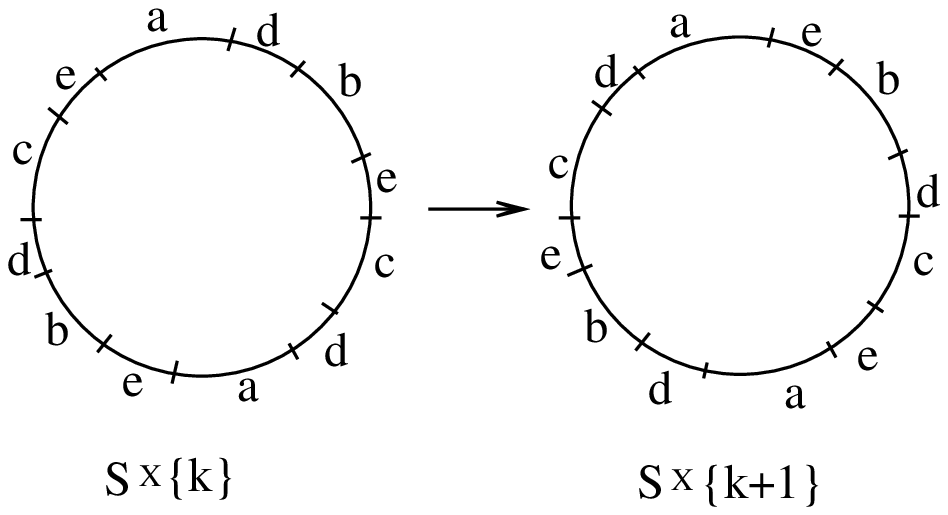}\\
{\it Переход от $k$-ой компоненты к $(k+1)$-ой: три красных
``двуугольника'' (символы $a,b$ и $c$) и два синих треугольника
(символы $d$ и $e$).}
\end{center}
}\label{RisGraphRauzy1}

{\bf Замечания.}

\begin{enumerate}
\item Из вида динамической системы следует, что при ``склейке'' всех
оттенков красного в один цвет, а всех оттенков синего -- в другой,
получится бинарное слово Штурма.

\item Поскольку все середины дуг одного цвета находятся в вершинах
правильного $m$-угольника, и оттенки внутри цвета также являются
вершинами правильных $k_i$-угольников, то по теореме о покрытии
множества целых арифметическими прогрессиями, найдутся два
одинаковых многоугольника.

\end{enumerate}

\subsection{Замечания о периодических сбалансированных словах над произвольным алфавитом}

В рассуждениях о непериодических сбалансированных словах мы
использовали всюдуплотность траектории для каждого слова $W_i$. В
периодическом случае эти рассуждения не проходят. Однако из
имеющейся динамической системы можно сконструировать периодическое
сбалансированное слово.

\begin{proposition}
Пусть динамическая система удовлетворяет условиям основной теоремы
с той лишь разницей, что $\alpha$ -- рациональное число. Тогда
соответствующее слово является сбалансированным периодическим.
\end{proposition}
\Proof Покажем сбалансированность по произвольному символу $a_i$.
В нашей динамической системе "склеим" все цвета, отличные от
цвета, соответствующего $a_i$. Нетрудно видеть, что в этом случае
система будет эквивалентна системе, подобной системе, порождающей
слова Штурма, но с рациональной величиной сдвига. По предложению
\ref{OsnTh} это слово является периодическим сбалансированным.
\Endproof

Возникает естествнный\par {\bf Вопрос:} {\it любое ли
периодическое слово может быть получено таким путем?}

Ответ на этот вопрос отрицательный. Действительно, как следует из
замечания к основной теореме, в случае алфавитов более, чем из
двух символов, для сбалансированных слов, порождающихя такой
динамической системой, обязательно существуют символы, имеющие
одинаковую плотность.

Но уже для алфавита из трех символов существует слово с попарно
различными плотностями символов:

$$
W=(abacaba)^\infty
$$

Несложно проверить, что $W$ является сбалансированным. Эту
конструкцию можно расширить на алфавиты из большего числа символов
рекуррентным образом:

\begin{equation}\label{Frankel}
W_k=(U_{k-1}a_kU_{k-1}), \ U_3=a_1a_2a_1a_3a_1a_2a_1
\end{equation}

В связи с этим сформулируем гипотезу Френкеля:

{\bf Гипотеза (Fraenkel).} {\it Единственное (с точностью до
перестановок) сбалансированное слово над алфавитом из $k\geq 3$
символов с попарно различными плотностями есть слово
(\ref{Frankel})}.

В настоящий момент гипотеза доказана для $k=3,4,5,6,7$.

\subsection{Сбалансированные слова и теорема Голода-Шафаревича.\label{GolodShaf}}

Е.С.Голод и И.Р. Шафаревич построили пример $k$-порожденных ненильпотентных алгебр таких, что каждая 
$(k-1)$-порожденная подалгебра нильпотентна. Нильпотентность всякой $1$-порожденной подалгебры  означает
условие нилевости. 

В ненулевых словах этих алгебр образующие должны быть перемешаны. Понятие сбалансированности может служить описанием процесса
перемешивания. 

Можно показать, что существует $2$-сбалансированное сверхслово $W$ с экспоненциальной функцией роста. Система определяющих соотношений
$u_i=0$, где $u_i$--не подслово $W$,   определяет мономиальую алгебру $A$ экспоненциального роста. Действуя, как и в работе Голода-Шафаревича \cite{GolShaf}
можно построить ненильпотентную ниль-алгебру, в которой все ненулевые слова $2$-сбалансированны. Вместе с тем, условие $1$-сбалансированности является чрезвычайно жестким. Наше описание сбаласнированных слов позволяет надеяться доказать следующую \par
{\bf гипотезу:}
{\it Ниль-алгебра, у которой все неприводимые слова $1$-сбалансированны, локально нильпотентна.}

%% file: chapter3.tex
\section{Слова с минимальной функцией роста.}\label{graph}
В этой главе изучается другое естественное обобщение слов Штурма --
слова с минимальной функцией роста, то есть с функцией роста,
удовлетворяющей соотношению $F_W(n+1)-F_W(n)=1$ при всех
достаточно больших $n$.

\subsection{Постановка задачи.}

Другим обобщением слов Штурма являются слова с минимальной
функцией сложности. Наша цель -- построить динамическую систему,
которая аналогичным образом могла бы порождать слова медленного
роста. В этой части мы будем рассматривать слова над произвольным
алфавитом $A=\{a_1,a_2,\ldots ,a_n\}$.

Пусть слово $W$ будет словом медленного роста, то есть, начиная с
некоторого $N$ верно $T_W(k+1)=T_W(k)+1, k\geq N$.

Рассмотрим $k$-графы этого слова для $k \geq N$. Так как ребра
графа соответствуют ($k+1$)-словам, то в этом графе $m$ вершин и
$m+1$ ребро и, соответственно, имеет  одну входящую и одну
выходящую развилку, которые, возможно, совпадают (если данная
вершина соответствует биспециальному слову).

Имеется два типа таких графов:

\vbox{
\begin{center}
\includegraphics{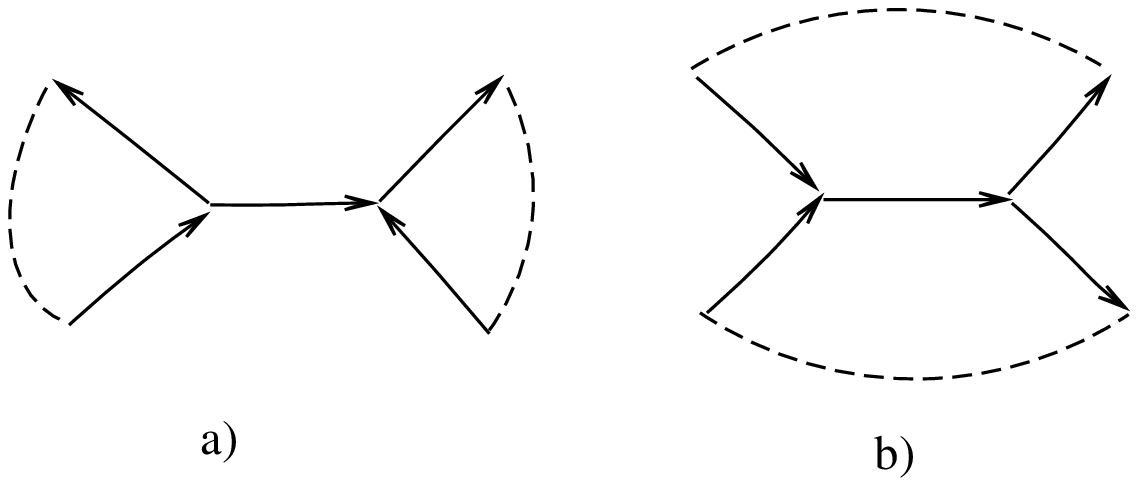}\\
{\it Два типа графов с одной входящей и одной входящей развилкой.}
\end{center}
}\label{RisGraphRauzy1}

Легко видеть, что граф типа a) не является сильно связным, а в
случае b) -- является. Слова с $k$-графами типа b) имеют вид
$u^{\infty}Vv^{\infty}$  и не являются равномерно-рекуррентными

Нас интересует случай a).

Назовем путь, ведущий из входящей развилки в выходящую {\it
перегородкой}, а два пути из выходящей развилки будем называть
{\it дугами}.

\begin{proposition}
Пусть $k$-граф слова $W$ имеет перегородку длины $l\geq 1$ и дуги
длины $r,s$. Тогда последователь $\Fol(G)$ имеет перегородку длины
$l-1$, дуги длины $r+1$ и $s+1$ и совпадает с $(k+1)$-графом слова
$W$.

\end{proposition}

Рассмотрим теперь предельный случай, когда перегородка
вырождается, то есть входящая развилка совпадает с выходящей. В
этом случае развилка соответствует биспециальному слову.

\begin{proposition}
Последователь графа с вырожденной перегородкой имеет две входящие
и две выходящие развилки

\end{proposition}

Пусть $u$ -- биспециальное подслово $W$, то есть при некоторых
$a_i,a_j,a_r,a_s\in A$ $a_i u,a_j u, ua_r, ua_s$ -- тоже подслова
$W$. Тогда развилками (входящими и выходящими) в $\Fol(G)$ будут
вершины соответствующие этим словам.

Таким  образом, для $(k+1)$-графа слова $W$ имеется $4$
возможности для удаления одного ребра из $\Fol(G)$,
соответствующего минимальному запрещенному $(k+2)$-слову: $a_i u
a_r, a_i u a_s, a_j u a_s, a_j u a_r$, В двух случаях мы получаем
сильно связный граф, а в двух -- не сильно связный.

Непосредственной проверкой доказывается
\begin{proposition}
Пусть в графе $G$ перегородка вырождена, а дуги имеют длину $r,s$
соответственно. Тогда $(k+1)$-граф имеет вид $(s-1,1,r+1)$ или
$(r-1,1,s+1)$.
\end{proposition}

\begin{proposition}\label{graph:3}
Любой граф $G$ с одной входящей и одной выходящей развилкой имеет
предшественника, причем только одного.
\end{proposition}

\Proof Пусть граф имеет вид $(l,r,s)$. Заметим, что графа с
$r=s=1$ не бывает. В случае, если $r,s>1$ граф предшественника
имеет вид $(l+1,r-1,s-1)$, если имеет вид $(l,1,s),s>1$ , то
предшественник -- $(0,l+1,s-1)$, если $(l,r,1),r>1$ -- то
$(0,s+1,l-1)$. Граф вида $(0,1,k)$ имеет предшественника
$(0,1,k-1)$.\Endproof

Ясно, что любой граф может быть получен из графа типа $(0,1,2)$
или $(0,2,1)$. Из предложения \ref{graph:3} следует, что для слова
с минимальной функцией роста $W$ существует такое слово Штурма $V$
и такие натуральные $n$ и $l$, что  $k$-графы для слова $W$
совпадают с $(k+l)$-графами слова Штурма $V$ для $k\geq n$.

\subsection{Свойства слов с минимальной функцией роста}

\begin{proposition}
Пусть $W$ -- слово с минимальной функцией роста, тогда существует
натуральное $k$, что для любого подслова $v\subset W$ такого, что
$|v|\geq k$, сущствует ровно два возвращаемых слова.
\end{proposition}

\Proof Рассмотрим слово Штурма, такое что эволюция $k$-графов Рози
слова $W$ и графов слова $V$ совпадают, начиная с некоторого $n$.
Тогда $k$-словам слова $W$ биективно соответствуют $(k+l)$-словам
слова $V$. Поскольку для подслов слова Штурма существует ровно $2$
возвращаемых слова, то это же верно и для подслов слова
$W$.\Endproof

Из этих же соображений доказывается
\begin{proposition}
рекуррентное слово с минимальной функцией роста равномерно
рекуррентно.
\end{proposition}
\Noproof

\subsection{Эволюции k-графов слов медленного роста}

В этой части мы рассматриваем эволюцию $k$-графов бесконечных слов
с минимальным ростом, то есть таких слов W, у которых функция
сложности растет минимально: $T_W(n+1)-T_W(n)=1$, начиная с
какого--то $n$.

В случае, когда $T_W(n+1)-T_W(n)=1$ для всех $n>0$, мы получим
эволюцию $k$-графов слов Штурма.

Так как ребра графа соответствуют ($k+1$)--словам, то в таком графе
$m$ вершин и $m+1$ ребро и, соответственно, имеет одну входящую и
одну выходящую развилку, которые, возможно, совпадают (если данная
вершина соответствует биспециальному слову).

Имеется два типа таких графов: в одном случае граф является сильно
связным, в другом -- нет.

Слова с $k$-графами типа b) имеют вид $u^{\infty}Vv^{\infty}$. Для
нас интересным является случай а).

Назовем путь, ведущий из входящей развилки в выходящую {\it
перегородкой}, а два пути из выходящей развилки будем называть
{\it дугами}.

\begin{proposition}
Пусть k-граф слова $W$ имеет перегородку длины $l\geq 1$ и дуги
длины $r,s$. Тогда последователь $\Fol(G)$ имеет перегородку длины
$l-1$, дуги длины $r+1$ и $s+1$ и совпадает с $(k+1)$-графом слова
$W$.
\end{proposition}
\Proof Непосредственно проверяется.\Endproof

\vbox{
\begin{center}
\includegraphics{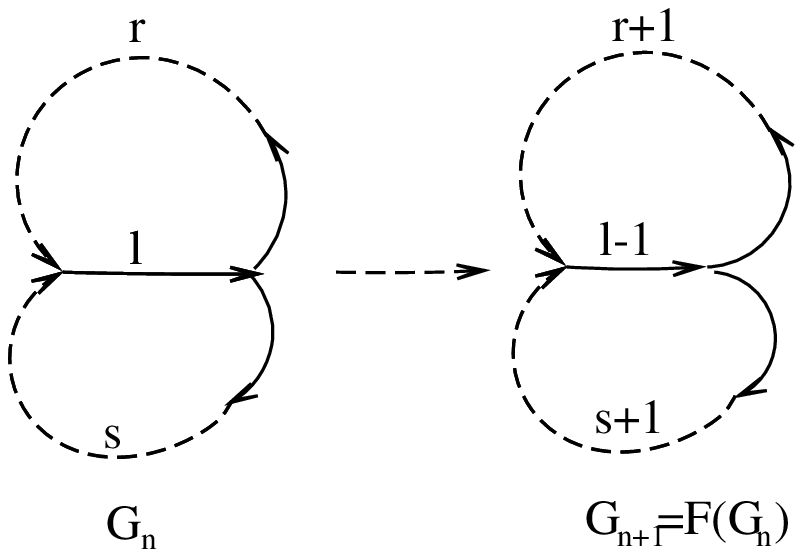}\\
{\it Переход от $G_n$ к $G_{n+1}$ для невырожденной развилки.}
\end{center}
}\label{RisGraphRauzy1}

Рассмотрим теперь предельный случай, когда перегородка
вырождается, то есть входящая развилка совпадает с выходящей. В
этом случае развилка соответствует биспециальному слову.

\begin{proposition}
Последователь графа с вырожденной перегородкой имеет 2 входящие и
две выходящие развилки

\end{proposition}

Пусть $u$ -- биспециальное подслово $W$, то есть при некоторых
$a_i,a_j,a_r,a_s\in A$ $a_i u,a_j u, ua_r, ua_s$ -- тоже подслова
$W$. Тогда развилками (входящими и выходящими) в $\Fol(G)$ будут
вершины соответствующие этим словам.

Таким  образом, для $(k+1)$--графа слова $W$ имеется четыре
возможности для удаления одного ребра из $\Fol(G)$,
соответствующего минимальному запрещенному $(k+2)$-слову: $a_i u
a_r, a_i u a_s, a_j u a_s, a_j u a_r$, В двух случаях мы получаем
сильно связный граф, а в двух -- не сильно связный.

\begin{proposition}

Пусть в графе $G$ перегородка вырождена, а дуги имеют длину $r,s$
соответственно. Тогда $(k+1)$-граф имеет вид $(s-1,1,r+1)$ или
$(r-1,1,s+1)$.

\end{proposition}

\begin{proposition}
Любой граф $G$ с одной входящей и одной выходящей развилкой имеет
предшественника, причем только одного.
\end{proposition}

\Proof Пусть граф имеет вид $(l,r,s)$. Заметим, что графа с
$r=s=1$ не бывает. В случае, если $r,s>1$ граф предшественника
имеет вид $(l+1,r-1,s-1)$, если имеет вид $(l,1,s),s>1$ , то
предшественник есть $(0,l+1,s-1)$, если $(l,r,1),r>1$ -- то
$(0,s+1,l-1)$.

Граф вида $(0,1,k)$ имеет предшественника $(0,1,k-1)$.

 Ясно, что любой граф может быть получен из графа типа
$(0,1,2)$ или $(0,2,1)$.

\subsection{Частоты подслов в словах медленного роста}
Напомним, что {\it частотой} какого-либо подслова $w\subset W$
бесконечного слова $W$ называется следующий предел (если он
существует):
\begin{equation}\label{dens}
\mu(w)=\lim_{n\to \infty} \mu_n(w),
\end{equation}
где $\mu_n(w)$ -- количество вхождений подслова $w$ в префикс $W$
длины $n$.

В случае слов Штурма предел существует всегда и равен длине
характеристического интервала для данного подслова, что следует из
равномерной распределенности последовательностей вида $\{\alpha
n+\beta\}$ для иррациональных $\alpha$. Докажем полезное

\begin{proposition}\label{equiv_freq}
Пусть $W$ -- слово с минимальной функцией роста и пусть $(v,a_i,u)$
-- ребро в его $k$--графе, где $k\geq N$. Тогда:

\begin{enumerate}

\item Любое конечное подслово слова $W$ обладает частотой, то есть
существует предел (\ref{dens}).

\item Если $v$ не правое специальное подслово, а $u$ -- не левое
специальное, то $\mu(u)=\mu(v)$.

\end{enumerate}
\end{proposition}

\Proof Из определения $k$ графов следует, что существует символ
$a_j\in A$ такой, что $ua_i=a_jv\in F_{n+1}(W)$. Так как $ua_k,
a_lv\notin F_W(n+1)$ для любых $a_k\neq a_i$ и $a_l\neq a_j$, то
$\mu(u)=\mu(ua_i)=\mu(a_jv)=\mu(v)$.

\Endproof

Следующее утверждение является обобщением т.н. {\it теоремы о трех
расстояниях } (three distance theorem).

\begin{theorem}
Пусть $W$ -- слово медленного роста над произвольным алфавитом $A$
с функцией роста $T_W(n)$ как в предложении \ref{equiv_freq}.
Тогда для любого $n\geq N$ подслова длины $n$ могут принимать не
более трех значений. Если они принимают три значения, то одно из
них является суммой двух других.
\end{theorem}

Эволюция графов с одной развилкой может быть обобщена на случай
слов с функцией сложности $T_W(n)$, удовлетворяющей соотношению
$$\lim_{n \to \infty} {T_W(n)/n}=1$$.

\subsection{Конструкция динамической системы}

Теперь мы покажем, как по соответствующему графу слов построить
динамическую систему, которая порождала бы данное слово.

Построение динамической системы мы осуществим двумя способами.
Первый способ опирается на уже известный результат о словах
Штурма. Второй способ -- более конструктивный и обобщающийся на
случай нескольких развилок.

Итак, пусть слово $W$ имеет минимальный рост, т.е., с
некоторого $k$ его $k$-графы слова $W$ имеют вид a).

Поскольку каждый граф имеет единственного предшественника, то
существует последовательность графов типа a) $G'_1, G'_2, \ldots,
G'_n$ таких, что $G'_{l+1}$ является предшественником $G'_l$, а
граф $G'_n$ совпадает с $k$-графом слова $W$.

Иными словами, существует эволюция $k$-графов типа a), начиная с
$k=1$, что $n$-ый граф в эволюции совпадает с $k$-графом слова
$W$, ($n+1$)-ый  совпадает с ($k+1$)-графом слова $W$ и т.д.

Последовательность графов такого типа однозначно соответствует
некоторому слову Штурма $V$.

Значит, существует взаимнооднозначное соответствие между
$k$-словами слова $W$  и $n$-словами слова Штурма $V$, при этом
соответствие продолжается на ($k+1$)-слова $W$ и ($n+1$)-слова $V$
и т.д.

Разберем сначала случай, когда $k=1$, то есть уже $1$-граф слова
$W$ имеет вид а).  Символам алфавита $A$ взаимнооднозначно
соответствуют $n$-слова некоторого слова Штурма $V$. По теореме
эквивалентности, слово $V$ порождается сдвигом окружности
$T_\alpha$. Как показано в части \ref{Chapter1}, $n$ - словам в
динамике соответствуют интервалы разбиения: слову $w=w_1w_2\ldots
w_n$ ($w_i\in \{a,b\}$) соответствует интервал
$$
I_w=T^{(-n+1)}_\alpha(I_{w_1})\cup T^{(-n+2)}_\alpha(I_{w_2}) \cup
\ldots \cup I_{w_n},$$ где $I_{w_i}$ -- характеристичесий интервал
для $w_i$.

Для сдвига окружности интервалы $I_w$, соответствующие $n$-словам,
будут иметь вид $I_w=(n_i\alpha,n_{i+1}\alpha)$, где $n_i$ --
некоторые целые числа.

Тогда характеристическим множеством для каждого символа из
алфавита $A$ будет являться интервал $I_w$, такой, что слово $w$
соответствует данному символу.  Ясно тогда, что слово $W$ будет
порождаться тем же сдвигом $T_\alpha$ и характеристическими
множествами $I_w$.

Теперь разберем общий случай. Пусть $k$-словам слова $W$
соответствуют $n$-слова слова Штурма $V$.  Точно также мы можем
построить соответствие между интервалами разбиения для $n$-слов
слова Штурма $I_w$ и $k$-словами слова $V$.  Построим
характеристическое множество для каждого символа из $A$.

А, именно, символу $a_i$ поставим в соответствие все интервалы,
которые соответствуют словам, начинающимся на $a_i$.  В этом случае
характеристическое множество для произвольного символа может быть
несвязным и представлять собой объединение нескольких интервалов.

Покажем, что при таком выборе характеристических множеств
найдется точка, чья эволюция будет совпадать с $W$.

Действительно, пусть слово $W$ начинается с некоторого $k$-слова
$w=w_1w_2\cdots w_k$. Ему мы поставим в соответствие интервал $I_w$.
Рассмотрим следующий $k+1$ символ $w_{k+1}$; этому ($k+1$)-слову мы
поставим в соответствие интервал $I_w\cup I_{w'}$, где $w'=w_2w_3
\cdots w_{k+1}$, и т.д.

Бесконечное пересечение интервалов, соответствующих префиксам
данного слова будет содержать точку, эволюция которой совпадает с
$W$.

Таким образом, доказана следующая

\begin{theorem}\label{MinGrow}
Пусть $W$ -- рекуррентное слово над произвольным конечным
алфавитом $A$. Тогда следующие условия на слово $W$ эквивалентны:
\begin{enumerate}

\item Существует такое натуральное $N$, что функция сложности
слова $W$ равна $T_W(n)=n+K$, для $n\geq N$  и некоторого
постоянного натурального $K$.

\item Существуют такое иррациональное $\alpha$ и целые $n_1,n_2,
\ldots,n_m$, что слово $W$ порождается динамической системой
$$
(\mathbb{S}^1,T_\alpha, I_{a_1},I_{a_2},\ldots ,I_{a_n},x),
$$
где
$T\alpha$ -- сдвиг окружности на иррациональную величину $\alpha$,
$I_{a_i}$ -- объединение дуг вида $(n_j\alpha,n_{j+1}\alpha)$.

\end{enumerate}

\end{theorem}

\subsection{Нормальные базисы граничных аллгебр.}

Для произвольной алгебры $A$ через $V_A(n)$ обозначим размерность
вектрного пространства, порожденного мономами длины не больше $n$.
Пусть $T_A(n)=V_A(n)-V_A(n-1)$. Если алгебра однородна, то $T_A(n)$
есть размерность векторного пространства, порожденного мономами
длины ровно $n$.

Известно (см. \cite{BBL}) что либо $\lim_{n\to
\infty}{(T_A(n)-n)}=-\infty$ (в этом случае есть альтернатива ({\it
Bergman Gap Theorem}): либо $\lim{V_A(n)}=C<\infty$ и тогда $\dim A
< \infty$, либо $V_A(n)=O(n)$ и алгебра имеет {\it медленный рост}),
либо  $T_A(n)-n < \Const$, либо, наконец, $\lim_{n\to
\infty}{(T_A(n)-n)}=\infty$.

В последнем случае для любой функции $\phi(n)\to \infty$ и любой
$\psi(n)=e^{o(n)}$ существует алгебра $A$ такая, что для
бесконечного множества натруальных чисел $n\in L\subset {\mathbb N}\
T_A(n)>\psi(n)$ и для бесконечного множества натуральных чисел $n\in
M\subset {\mathbb N}\ T_A(n)<n+\phi(n)$ (см. \cite{BBL}).

Таким образом, в случае $\lim_{n\to \infty}{(T_A(n)-n)}=\infty$ рост
может быть хаотичным и мы этот случай не рассматриваем. Случай,
когда $T_A(n)-n < \Const$ (т.е. случай алгебр {\it медленного
роста}) исследовался Дж.~Бергманом и Л.~Смоллом. Нормальные базисы
для таких алгебр исследованы в работе \cite{BBL}. Назовем алгебру
{\it граничной}, если $T_A(n)-n<\Const$. Наша цель состоит в
описании нормальных базисов граничных алгебр.

Прежде всего отметим, что такое описание сводится к мономиальному
случаю. В самом деле, пусть $a_1,\ldots,a_s$ -- образующие алгебры
$A$. Порядок $a_1\prec\ldots\prec a_s$ индуцирует порядок на
множестве слов алгебры $A$ (сперва по длине, затем
лексикографически). Назовем слово {\it неуменьшаемым} (или {\it
неприводимым}), если его нельзя представить в виде линейной
комбинации меньших слов. Множество неуменьшаемых слов обазует {\it
нормальный базис} алгебры $A$ как векторного пространства.
Рассмотрим фактор свободной алгебры $k<\hat{a}_1,\ldots,\hat{a}_s>$
($k$~-- основное поле) по множеству слов из норомального базиса.
Получится мономиальная алгебра с тем же нормальным базисом и, стало
быть, с той же функцией роста.

Назовем {\it обструкцией} приводимое (то есть уменьшаемое) слово $u$
такое, что любое его подслово неприводимо. Сверхсловом в алгебре $A$
(правым, левым, двусторонним) называется сверхслово $W$ такое, что
любое его конечное подслово ненулевое. Аналогично определяется {\it
неприводимое сверхслово} в алгебре $A$.

\begin{theorem}[Описание нормальных базисов граничных алгебр]\label{ThMainLast}
Пусть $A$ -- граничная алгебра, $a_1,\ldots, a_s$ -- ее образующие.
Тогда имеют место два случая, каждый из которых имеет свое описание.

{\bf Случай 1.} Алгебра $A$ не содержит равномерно-рекуррентного
непериодического сверхслова. В этом случае нормальный базис алгебры
$A$ состоит из множества подслов следующего множества слов:

\begin{enumerate}
\item Одно слово вида $W=u^{\infty/2}cv^{\infty/2} \neq u^{\infty}$

\item Произвольный конечный набор $\mu$ конечных слов

\item Множество слов вида $u_i^{\infty/2}c_i$,\ $i=1,\ldots, r_1$

\item Множество слов вида $d_iv_i^{\infty/2}$,\ $i=1,\ldots, r_2$

\item Множество слов вида $e_j{(R_j)}^{k} f_j$,\ $k\in \mathbb{K}_j \subseteq
\mathbb{N}$, $j=1,\ldots,r_3$

\item Множество слов вида $W_{\alpha} = E_{\alpha}u^{n_\alpha}cv^{m_\alpha}
F_\alpha$. При этом существуют такое $c>0$, что для любого $k$
количество слов $W_\alpha$ длины $k$ меньше $c$.
\end{enumerate}

{\bf Случай 2.} Алгебра $A$ содержит равномерно-рекуррентное
непериодическое сверхслово $W$. В этом случае нормальный базис
алгебры $A$ состоит из множества подслов следующего семейства слов,
включающее в себя:

\begin{enumerate}
\item Некоторое равномерно рекуррентное слово $W$ с функцией роста
$T_W(n)=n+const$ для всех достаточно больших $n$. Описание таких
слов дано в теореме \ref{MinGrow}.

\item Произвольный конечный набор $\mu$ конечных слов

\item Множество слов вида $u_i^{\infty/2}c_i$, $i=1,\ldots,s_1$

\item Множество слов вида $d_iv_i^{\infty/2}$, $i=1,\ldots,s_2$

\item Множество слов вида $e_j(R_j)^{k}f_j$, $k\in \mathbb{K}_j \subseteq
\mathbb{N}$,\ $j=1,\ldots,s_3$

\item Множество слов вида $L_iO_iW_i$,\ $i=1,\ldots,k_1$. При этом
$W_i$ -- сверхслово, эквивалентное $W$ и $O_iW_i$ имеют вхождение
только одной обструкции (а, именно, $O_i$).

\item Множество слов вида $W'_jO'_jL'_j$, $j=1,\ldots,k_2$. При этом
$W'_j$ -- сверхслово, эквивалентное $W$ и $W'_jO'_j$ имеют вхождение
только одной обструкции (а, именно, $O'_j$).

\item Конечное множество серий вида: ${h^1}_iT_i{h^2}_i$,
$i=1,\ldots,s$. При этом:

a) слово $T_i$ содержит вхождение ровно двух обструкций
$O_{i_1},O_{i_2}$ (возможно, перекрывающихся)

b) Для некоторого $c>0$ $|h^1_i|+|h^2_i|<c$ при всех $i$.

с) Существует $m>0$ такое, что для любого $k$ имеется не более $k$
подслов длины $m$, вида $(8)$ и не являющихся подсловами слова $W$.

\end{enumerate}
\end{theorem}

\Proof Алгебру $A$ считаем мономиальной. Назовем слово $v$ алгебры
$A$ {\it хорошим}, если для любого $n$ существуют сколь угодно
длинные слова $w_1,w_2$, $|w_1|>n,|w_2|>n$ такие, что $w_1vw_2$ есть
подслово алгебры $A$. Обозначим $T_{RL}(n)$  количество хороших слов
длины $n$. Известно, что если $T_{RL}(n)=T_{RL}(n+1)$ при некотором
$n$, то алгебра $A$ имеет медленный рост (см. \cite{BBL}).  В силу
граничности алгебры $A$ $T_{RL}(n)\geq T_{RL}(n+1)+1$,  при всех
достаточно больших $n$ неравенство превращается в равенство (иначе
$\lim_{n\to \infty} (T_{RL}(n)-n)= \infty$ и, как следствие
$\lim_{n\to \infty} (T_{A}(n)-n)= \infty$).

При этом граф Рози имеет развилку и, как следствие, два цикла,
эволюция графа Рози устроена следующим образом: либо граф теряет
сильную связность и имеет вид a). В этом случае его дальнейшая
эволюция однозначна, она отвечает слову вида $u^{\infty/2} c
v^{\infty/2}$ и мы имеем случай $1$, либо граф Рози все время
остается сильно связным. Тогда эволюция связной компоненты, в
которой есть развилка, асимптотически эквивалентна эволюции графа
Рози некоторого слова Штурма. Если оно имеет вид  $u^{\infty/2} c
u^{\infty/2}$, то имеет место случай $1$, иначе оно
равномерно-рекуррентно и имеет место случай $2$.

В случае 1 все обструкции для слова $u^{\infty/2} c u^{\infty/2}$
имеют ограниченную длину (см. \cite{BBL}). Можно сделать следующее

\medskip
{\bf Наблюдение.}\ {\it Количество ненулевых слов длины $n$, не
являющихся подсловами слова $u^{\infty/2}cu^{\infty/2}$ не
превосходит константы (не зависящей от $n$).}
\medskip

Напомним предложение из работы \cite{BBL}:

\begin{proposition}
Пусть $SW=WT$. Тогда $W$ имеет вид: $s^k s_1$, где $s_1$ -- начало
слова $s$.
\end{proposition}

Из данного предложения и только что сделанного наблюдения вытекает

\begin{corollary}
Пусть $|u|=k$, $|v_1|=|v_2|=l$. Тогда либо количество подслов длины
$k+m$, (где $m\le k$) слова $v_1uv_2$ не менее $m+1$, либо
$u=s^ks'$, для некоторого слова $s$, при этом $s'$ -- начало $s$ и
$v_1=v'_1s'$.
\end{corollary}

Из данного следствия и наблюдения получается

\begin{proposition}
Существует константа $K$, зависящая только от граничной алгебры $A$,
такая, что для любой обструкции $O$ в слове $W$ либо при некотором
$m$ для любых $v_1,v_2$, $|v_1|\geq m$, $|v_2|\geq m$ $v_1Ov_2$
является нулевым словом алгебры (число $m$ не зависит от выбора
обструкции), либо $|v|\leq K$.
\end{proposition}

Из данного предложения следует, что словами алгебры $A$ являются
либо слова, содержащие не более двух обструкций, причем каждая
обструкция находится на ограниченном расстоянии от одного из концов,
либо подслова слов вида $R_iu^k_iT_i$. А все такие типы слов описаны
в условии теоремы \ref{ThMainLast}. \Endproof

%% file: chapter4.tex
\section{Перекладывания отрезков и символическая динамика.}
\subsection{Введение и постановка задачи.}
Данная глава посвящена изучению слов, получаемых из перекладывания
отрезков. Каждое такое слово $W$ удовлетворяет сотношению:
$$
F_W(n+1)-F_W(n)=k, \  k\geq N,
$$

которое означает, что количество $n$-подслов для достаточно больших
$n$ увеличивается не постоянную величину $k$. Мы ограничиваемся
случаем, когда $W$ рекуррентно. Наша цель -- получить комбинаторный
критерий порождаемости слова $W$ перекладыванием отрезков.

Напомним (см. \ref{PerOtr}), что если перекладывание $k$ отрезков
является регулярным, то эволюция любой точки является словом с
линейной функцией сложности: $T_{U(x)}(n)=n(k-1)+1$.

\subsection{Необходимые условия для порождаемости слова перекладыванием отрезков}
Пусть слово $W$ является эволюцией точки $x\in [0,1]$ при
перекладывании $k$ отрезков и характеристичеких множествах
$U_1,U_2,\ldots,U_n$. Каждое характеристическое множество $U_i$
является объединением нескольких непересекающихся интервалов или
полуинтервалов.

Как было показано в части \ref{Sootv:2},  подслова длины $k$
взаимнооднозначно соответствуют $k$-разбиениям характеристических
множеств. Поскольку точка границы одномерного множества может
являться границей только для двух множеств, то  $k$-подслово слова
$W$ может иметь максимум два продолжения. Мы получили первое
необходимое условие
 для того, чтобы слово порождалось перекладыванием отрезков:

\begin{proposition}
Пусть слово $W$ порождается перекладыванием отрезков. Тогда для
некоторого $N$ все специальные
 слова подслова длины, не меньше $N$ должны иметь валентность, равную $2$.
\end{proposition}

Это условие аналогично тому, что начиная с какого-то $N$ все
$k$-графы слова $W$ ($k\geq N$) должны иметь входящие и исходящие
развилки степени $2$.

Следующие условия мы получим на языке графов Рози. Пусть
рекуррентное слово $W$ имеет функцию роста $F_W(n)=Kn+L$ для $n>N$ и
порождается перекладыванием отрезка. Рассмотрим эволюцию $k$-графов
Рози слова $W$, начиная с $k\geq N+1$. Как было показано выше, все
входящие и выходящие развилки имеют степень $2$, поэтому у всех
$k$-графов есть ровно $K$ входящих и $K$ исходящих развилок.

В минимальном случае, описанном в предыдущей главе, когда
$F_W(n)=n+L$, графы Рози имеют ровно одну входящую и одну исходящую
развилку. Если входящая и выходящая развилка совпадают, то в этом
случае выбор удаляемого ребра в последователе $D(G)$ определяется
однозначно из условия сильной связанности графа.

Более интересным случаем являются графы слов, в которых есть более
одной развилки. Рассмотрим его подробно. При переходе от графа $G_n$
к $G_{n+1}$ возможны следующие варианты:

\begin{enumerate}
\item В графе $G_n$ нет сцепленных циклов (то есть, нет входящих
развилок, являющихся одновременно исходящими). В этом случае граф
$G_{n+1}$ совпадает с последователем $D(G_n)$.

\item В графе $G_n$ одна развилка является одновременно входящей и
выходящей. Граф последователя $D(G_n)$ в этом случае имеет три
развилки, так как одна развилка размножилась. Следовательно, граф
$G_{n+1}$ получается из последователя $D(G_n)$ путем удаления одного
ребра, соответствующего минимальному невстречающемуся слову.

\item В графе есть две развилки или более развилок, являющихся одновременно входящими
и выходящими. Граф $G_{n+1}$ получается из $D(G_{n})$  удаления двух
или более ребер, соответствующих минимальным невстречающимся словам.

\end{enumerate}

Поскольку слово $W$ рекуррентно, то из предложения \ref{reccur}
следует, что при удалении ребер должна сохраняться сильная
связность, то есть из любой вершины можно по стрелкам перейти в
любую другую.

Рассмотрим подробнее второй случай. Пусть в $G_k$ есть одна двойная
развилка. Это означает, что в $W$ есть ровно одно биспециальное
подслово $w$ длины $k$. Значит, существуют такие $a_i,a_j,a_k,
a_l\in A$, что $a_iw,a_jw,wa_k,wa_l$ -- подслова $W$. Тогда
$(k+1)$-граф $G_{k+1}$ получается из последователя путем удаления
какого-то ребра, соответствующего одному из четырех слов:
$a_iwa_k,a_iwa_l, a_jwa_k, a_jwa_l$. Рассмотрим интервал, являющийся
характеристическим множеством для слова $I_w=[x_w,y_w]$.

Так как $w$ -- правое специальное слово, то $I_w\subset
T^{-1}(I_{a_k}\cup I_{a_l})$, так как $w$ -- левое специальное, то
$I_w\subset T(I_{a_i}\cup I_{a_j})$.

Пусть точка $A\in [0,1]$ делит $I_w$ на два интервала, образы
которых лежат в  $I_{a_k}$  и $I_{a_l}$ соответственно, а точка
$B\in [0,1]$ -- делит на интервалы, прообразы которых лежат в
$I_{a_i}$ и $I_{a_j}$ соответственно.

Выбор минимального невстречающегося  слова, а , значит, удаляемого
ребра, определяется взаиморасположением точек $A$ и $B$, а также
сохранением или сменой ориентации отображения на этих множествах.
Итого, имеется 8 вариантов, которые разбиваются на четыре пары,
соответствующие одинаковым наборам слов:

\begin{enumerate}

\item $B<A$, $T^{-1}([x_w,B])\subset I_{a_i}$, $T([x_w,A]\subset
I_{a_k} )$

\item $B<A$, $T^{-1}([x_w,B])\subset I_{a_j}$, $T([x_w,A]\subset
I_{a_k} )$

\item $B<A$, $T^{-1}([x_w,B])\subset I_{a_i}$, $T([x_w,A]\subset
I_{a_l} )$

\item $B<A$, $T^{-1}([x_w,B])\subset I_{a_j}$, $T([x_w,A]\subset
I_{a_l} )$

\item $B>A$, $T^{-1}([x_w,B])\subset I_{a_i}$, $T([x_w,A]\subset
I_{a_k} )$

\item $B>A$, $T^{-1}([x_w,B])\subset I_{a_j}$, $T([x_w,A]\subset
I_{a_k} )$

\item $B>A$, $T^{-1}([x_w,B])\subset I_{a_i}$, $T([x_w,A]\subset
I_{a_l} )$

\item $B>A$, $T^{-1}([x_w,B])\subset I_{a_j}$, $T([x_w,A]\subset
I_{a_l} )$

\end{enumerate}

\vbox{
\begin{center}
\includegraphics{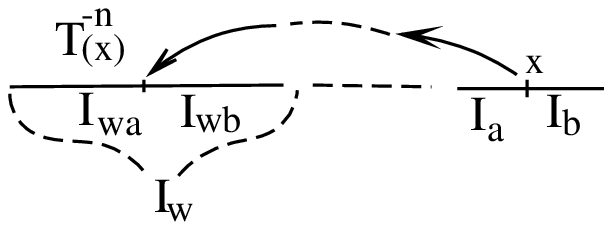}\\
{\it Разбиение характеристического множества для слова $w$ на $n$-м
шаге.}
\end{center}
}\label{RisGraphRauzy1}

Варианты $1$ и $5$ соответствуют запрещенному слову $a_iwa_k$.

Варианты $2$ и $6$ соответствуют запрещенному слову $a_iwa_l$.

Варианты $3$ и $7$ запрещенному слову соответствуют $a_jwa_k$.

Варианты $4$ и $8$ соответствуют запрещенному слову $a_jwa_l$.

Две пары вариантов соответствуют одновременной смене или сохранению
ориентации отображения на характеристических интервалах, а две пары
-- различным ориентациям.

В случае, когда перекладывание отрезков  сохраняет ориентацию, у нас
остается только два возможных варианта для удаляемого ребра.

В случае, когда мы разрешаем отобажению менять ориентацию отрезков
при перекладывании, возможны все четыре варианта.

Введем понятие {\it размеченного графа Рози}. Граф Рози называется
{\it размеченным}, если

\begin{enumerate}
\item Ребра каждой развилки помечены символами $l$ (``left'') и $r$ (``right'')

\item Некоторые вершины помечены символом ``--''.
\end{enumerate}

{\it Последователем} размеченного графа Рози назовем ориентированный
граф, являющийся его последователем как графа Рози, разметка ребер
которого определяется по правилу:
\begin{enumerate}

\item Ребра, входящие в развилку должны быть помечены
теми же символами, как и ребра, входящие в любого левого потомка
этой вершины;

\item Ребра, выходящие из развилки должны быть помечены
теми же символами, как и ребра, выходящие из любого правого потомка
этой вершины;

\item Если вершина помечена знаком ``--'', то все ее правые
потомки также должны быть помечены знаком ``--''.

\end{enumerate}

{\bf Замечание.} Поясним смысл разметки графа. Пусть ребра входящей
развилки соответствуют $a_i$ и $a_j$, символы $l$ и $r$
соответствуют левому и правому множеству в паре
$(T(I_{a_i}),T(I_{a_j}))$. Если символы $a_k$  и $a_l$ соответствуют
ребрам исходящей развилки, то символы $l$ и $r$ ставятся в
соответствии с порядком ``лево-право'' в паре $(I_{a_k},I_{a_l})$.
Знак ``--'' ставится в вершине, если характеристическое множество,
ей соответствующее, принадлежит интервалу перекладывания, на котором
меняется ориентация.

Условие для перехода от графа $G_n$ к $G_{n+1}$:

\begin{proposition}
\begin{enumerate}
\item Если в графе нет двойных развилок, соответствующих
 биспециальным подсловам, то при переходе от $G_n$ к $G_{n+1}$ имеем $G_{n+1}=D(G_n)$;

\item Если вершина, соответствующая биспециальному слову не
 помечена знаком ``--'', то ребра, соответствующие запрещенным
 словам выбираются из пар $lr$ и $rl$
\item Если вершина помечена  знаком ``--'', то удаляемые ребра
должны выбираться из пары $ll$ или $rr$.
\end{enumerate}
\end{proposition}

Назовем эволюцию размеченных графов Рози {\it правильной}, если {\bf
правила 1} и {\bf 2} выполняются для всей цепочки эволюции графов,
начиная с $G_1$, назовем эволюцию {\it асимптотически правильной},
если {\bf правила 1} и {\bf 2} выполняются, начиная с некоторого
$G_n$. Будем говорить, что эволюция размеченных графов Рози {\it
ориентированна}, если в $k$-графах нет вершин, помеченных знаком
``--''.

С помощью введенных определений мы можем сформулировать следующие
условия на слово, порождающееся перекладыванием отрезков:

\begin{theorem}\label{Osn}
Пусть рекуррентное слово $W$
\begin{enumerate}
\item  порождается перекладыванием отрезков. Тогда это слово
обеспечивается асимптотически правильной эволюцией размеченных
графов Рози.
\item  порождается перекладыванием отрезков  с сохранением ориентации.
 Тогда это слово обеспечивается асимптотически правильной
 ориентированной эволюцией размеченных графов Рози.
\end{enumerate}
\end{theorem}

\subsection{Построение динамической системы}
Покажем, что условия теоремы \ref{Osn} являются достаточными для
того, чтобы слово порождалось перекладыванием отрезков. Сначала мы
покажем, что слово $W$, удовлетворяющее условиям теоремы \ref{Osn}
может являться эволюцией некоторой точки при кусочно-непрерывном
преобразовании отрезка $T:I\to I$ следующего вида:

\begin{enumerate}
\item $I=[x_0,x_1]\cup [x_1,x_2] \cup \ldots \cup [x_{n-1},x_n]$,
$x_0=0$, $x_n=1$.

\item $I=[y_0,y_1]\cup [y_1,y_2] \cup \ldots \cup [y_{n-1},y_n]$,
$y_0=0$, $y_n=1$.

\item $\sigma \in {S_n}$ -- некоторая перестановка из $n$
элементов.

\item $T$ отображает $(x_i,x_{i+1})$ на
$(y_{\sigma(i)},y_{\sigma(i)+1})$ непрерывно и взаимнооднозначно.

\end{enumerate}

Затем будет показано, как случай кусочно-непрерывного преобразования
 можно свести к перекладыванию отрезков.

Строить кусочно-непрерывное преобразование $T$ мы будем поэтапно.
При первой итерации мы разбиваем отрезок произвольным образом на
отрезки, которые соответствуют соответствующим символам. Для
построение отображения нам достаточно определить траекторию этих
точек и затем из соображений рекуррентности продолжить по
непрерывности на весь отрезок.

{\bf Замечание.} Из непрерывности и биективности следует, что
отображения на интервалах являются монотонными функциями. Мы
рассмотрим два случая:

\begin{enumerate}

\item Все отображения интервалов являются одновременно возрастающими
функциями. Назовем такое преобразование {\it сохраняющим
ориентацию}.

\item Во втором случае встречаются как возрастающие, так и убывающие
отображения интервалов. Этот случай назовем {\it не сохраняющим
ориентацию}.
\end{enumerate}

Разобьем отрезок $I=[0,1]$ на $n=\Card A$ интервалов произвольным
образом, которые будут являться характеристическими множествами для
символов алфавита $A$: $[0,1]=I_{a_1}\cup I_{a_2}\cup \ldots \cup
I_{a_n}$. Соответствие интервалов символам алфавита будет определено
в соответствии с предложением \ref{Sosed}

Будем считать, что отображение $T$ непрерывно на каждом множестве
$I_{a_i}$, то есть точками разрыва отображения $T$ могут быть только
конечные точки характеристических множеств. Интервалы
характеристических множеств или их образы, имеющие общую граничную
точку, будем называть {\it соседними}.

{\bf Замечание.} Мы всегда можем  расширить алфавит таким образом,
чтобы в новом алфавите характеристические множества были устроены
ровно так, как описано выше. Графы $k$-слов для старого алфавита
будут соответствовать $1$-графам нового алфавита, а эволюции графов,
начиная с этого момента, будут совпадать.

Непосредственно проверяется следующее
\begin{proposition}
Преобразования $Т, T^{(-1)}$ переводит два соседних множества  либо
в соседние множества, либо их образы не могут целиком покрывать
интервал.
\end{proposition}

\begin{proposition}\label{Sosed}
Интервалы разбиения могут быть поставлены в соответствие символам
алфавита таким образом, что если $w$ -- специальное правое $1$-слово
и $wa_i,wa_j$ -- подслова, то $I_{a_i}$  и $I_{a_j}$ -- соседние
множества. Аналогично, если $w$ -- специальное левое $1$-слово и
$a_kw,a_lw$ -- подслова, то $I_{a_k}$  и $I_{a_l}$ -- также соседние
множества.
\end{proposition}

\Proof Рассмотрим произвольную развилку в $1$-графе. Пусть из нее
стрелки ведут в вершины, соответствующие символам $a_i$ и $a_j$.
Тогда положим $I_{a_i}=[x_0,x_1]$, $I_{a_j}=[x_1,x_2]$,

Так как характеристические множества являются интервалами, то на них
можно естественным образом ввести отношение порядка (то есть
$I_{a_i}< I_{a_j}$, если $x_i<x_j$ ), точно также можно ввести
отношение порядка и на их образах. Если при преобразовании $T$ для
какой-то пары характеристических множеств порядок меняется, будем
говорить, что на этой паре преобразование меняет ориентацию.

Рассмотрим образы интервалов при отображении $T^{-1}$. Ясно, что
если символ $a_i$ не является правым специальным $1$-словом и за ним
однозначно следует символ $a_j$, то $I_{a_i}\subset
T^{-1}(I_{a_j})$, а если не является левым специальным и ему всегда
предшествует символ $a_k$, то $T^{-1}(I_{a_i})\subset I_{a_k}$.

В случае, если символ $a_i$ является правым специальным $1$-словом и
за ним могут идти символы $a_j$ и $a_k$, то $I_{a_i}\subset
T^{-1}(I_{a_j}\cup I_{a_k})$, а если левым специальным, то
$T^{-1}(I_{a_i})\subset I_{a_j}\cup I_{a_k}$.

Обозначим через $I^n$ множество образов концов интервалов при
отображении $T^{-n}$, $I^{0}$ -- множество концов интервалов
характеристических множеств, то есть $I^{0}=\{x_0,x_1,\ldots,
x_n\}$.

Как следует из рассуждений в части \ref{Sootv:2}, множество,
отвечающие слову $w=w_1w_2\cdots w_n$ есть $I_w=I_{w_n}\cap
T^{-1}(I_{w_{n-1}})\cap \ldots \cap T^{(-n+1)}(I_{w_1})$,
соответственно, множество граничных точек множеств, соответствующих
словам длины $n$ есть в точности $I^{0}\cup I^{1}\cup \ldots \cup
I^{(n-1)}$.

Если правое специальное слово не является биспециальным (то есть не
является одновременно и левым специальным), то положение точки,
которая делит данное характеристическое множество, несущественно и
его можно выбирать произвольно.

Для случая сохраняющего ориентацию преобразования удаления должно
соответствовать сохраняющим ориентацию правилам.

В случае не сохраняющего ориентацию преобразования нам просто
необходимо, чтобы внутри отображаемых отрезков при эволюции
происходило конечное число ``переломов''.

Таким образом мы  можем определить образы точек на некотором
подмножестве $N\subset I$. Из построения следует, что существуют
такие интервалы $I_k=(x_k,x_{k+1})$, внутри которых наше
преобразование монотонно. Мы всегда можем продолжить его по
непрерывности до отображения отрезка в себя. Построенное
кусочно-непрерывное преобразование и есть искомое. Обозначим его
$T$. Отметим, что начальная точка, эволюцией которой является
искомое слово $W=\{w_n\}$, получается как предельная
последовательность множества вложенных отрезков, соответствующих
префиксам
 $w_0, w_0w_1, w_0w_1w_2, \ldots$ и т.д.

Нами доказана

\begin{theorem}
Чтобы равномерно-рекуррентное слово $W$ порождалось
кусочно-непрерывным преобразованием отрезка с сохранением
ориентации, необходимо и достаточно чтобы слово обеспечивалось
асимптотически правильной ориентированной эволюцией размеченных
графов Рози.
\end{theorem}

В случае меняющего ориентацию преобразования получается

\begin{theorem}
Чтобы равномерно-рекуррентное слово $W$ порождалось
кусочно-непрерывным преобразованием отрезка с сохранением
ориентации, необходимо и достаточно чтобы слово обеспечивалось
асимптотически правильной эволюцией размеченных графов Рози.
\end{theorem}

\subsection{Эквивалентность множества р.р слов, порождаемых
 кусочно-непрерывным преобразованием,
множеству слов, порождаемых перекладыванием отрезков.}

Покажем вначале, как можно перейти к динамике, в которой почти все
точки (в смысле меры Лебега) имеют различные существенные эволюции.
Рассмотрим существенные эволюции точек при преобразовании $T$.

Рассмотрим кусочно-непрерыное преобразование отрезков. Из теоремы о
существовании инвариантной меры (см. \cite{Sin}) следует, что любое
отображение компакта имеет инвариантную вероятностную меру. Значит,
отображение $T$ имеет инвариантную меру $\mu$. Следовательно, на
отрезке мы можем ввести новую метрику $d(x_1,x_2)=\mu((x_1,x_2))$.
Отметим, что из того, что $\mu(x_1,x_2)>0$ не следует, что точки
имеют различную существенную эволюцию, поскольку построенное
отображение может быть разрывным. Кроме того, нам потребуется
подобрать подходящую меру.

Назовем разбиение отрезка на характеристические множества {\it
хорошим}, если для каждого символа его характеристическое множество
является выпуклым (то есть представляет собой отрезок, интервал или
поуинтервал) и если между двумя точками одноцветного множетсва есть
точка разрыва -- то их цвета их образов различны.

Для произвольного разбиения $U_1,U_2,\ldots,U_n$ назовем его {\it
доразбиением} любое разбиение $V_1, V_2,\ldots, V_m$  такое, что
каждое характеристическое множество $U_i$ является объединением
нескольких множеств $V_j$: $U_i=V_{i_1}\cup V_{i_2}\cup V_{i_j}$.
Эволюция доразбиения есть слово над алфавитом из большего числа букв
и эволюция разбиения получается из него путем склейки букв. Ясно,
что любое разбиение имеет хорошие доразбиения.

Пусть $W$ -- р.р. слово, соответствующее разбиению $\cal U$,
$\widehat{W}$ -- слово, соответствующее динамике с той же начальной
точкой и подразбиением ${\cal U}'$. Склеивающий морфизм алфавитов
индуцирует естественным образом морфизм слов $\pi: \widehat{W}\to
W$. $\widehat{W}$, вообще говоря, может не быть
равномерно-рекуррентным, но существует р.р. слово $\widehat{U}$
такое, что $\widehat{U}\preceq\widehat{W}$. Тогда
$U=\pi(\widehat{U})\preceq W=\pi(\widehat{W})$ и, значит, $U$
является р.р. словом (см. \ref{ThUnifRec}).

В соответствующем $\widehat{U}$ замыкании орбиты будут точки,
соответсвующие эволюции с проекцией $W$. Действительно, пусть $w$ --
произвольное подслово $W$, оно встречается в любом подслове $W$
длины $\ge k(w)$. Следовательно для любого подслова $\widehat{U}$
длины $\ge k(w)$ существует подслово $\hat{w}$ такое, что
$\pi(\hat{w})=w$. Следовательно, для любой близкой к нулевой позиции
$W$ существует $\widehat{V}\equiv \widehat{U}$ такое, что
$\pi(\widehat{V})$ будет также близко. Это позволяет использовать
соображения компактности. Соответсвующее слово $\widehat{W}'$ будет
р.р. словом.

Следующее предложение мы сформулируем для хороших разбиений.

\begin{proposition}
Пусть точки $x_1<x_2$ имеют одинаковые эволюции и преобразование $T$
непрерывно на интервале $(x_1,x_2)$. Тогда любая точка интервала
$(x_1,x_2)$ имеет такую же эволюцию.
\end{proposition}
\Proof Из того, что $x_1,x_2$ имеют одинаковую эволюцию следует, что
$T(x_1),T(x_2)$ тоже имеют одинаковую эволюцию и образом интервала
$(x_1,x_2)$ является интервал $(T(x_1),T(x_2))$.\Endproof

%Рассмотрим топологию, склеивающую точки с
%одинаковой существенной эволюцией и построим соответствующую
%фактор-динамику.  Нашем случае целые отрезки окажутся склеенными.
%Однако таких отрезков не более, чем счетное множество. В то же
%время, если $W$ -- непериодическое равномерно-рекуррентное слово, то
%количество эквивалентных (то есть имеющих одинаковый набор конечных
%подслов) равномерно-рекуррентных слов имеет мощность континуума (см.
%\cite{BBL}). Это означает, что существует несчетное множество точек
%вне этой системы отрезков и сингулярность меры на эволюцию этих
%точек влияние не оказывает. Таким образом, найдется эквивалентное
%слово $W'~W$, отвечающее перекладыванию отрезков. Но тогда в этом
%перекладывании встретятся {\it все} слова , эквивалентные $W'$, в
%том числе и слово $W$ (в смысле существенной эволюции).

Пусть $W$~-- непериодическое р.р. слово, наблюдающееся при кусочном
гомеоморфизме системы отрезков. Рассмотрим пространство
последовательностей с оператором сдвига. Для $W$ существует
минимальное замкнутое инвариантное множество $N_W$, отвечающее $W$
\cite{BBL}. Каждой точке $N_W$ отвечает множество точек системы
отрезков с данной существенной эволюцией. Далее, на пространстве
$N_W$ имеется вероятностная инвариантая относительно оператора
сдвига мера, которая индуцирует на системе отрезков естественным
образом вероятностную инвариантую меру (\cite{Sin}).

Если $I$ -- интервал, соответствующий слову $u$ и $N$ -- множество
точек замыкания орбиты $W$ начинающихся с подслова $u$, то
$\mu(I)=\nu(N)$. Если $I$ не пересекается с любым таким интервалом,
то $\mu(I)=0$.

Мера $\mu$ инвариантна и индуцирует полуметрику (полагаем длину
любого интервала равной его мере).

Рассмотрим топологию, склеивающую точки с нулевым расстоянием и
построим соответствующую фактор-динамику. Полученное отображение и
будет кусочной изометрией, т.е. перекладыванием отрезков.

Каждый такой склееный интервал содержится в некотором максимальном
склееном интервале, количество таких максимальных интервалов счетно
и их общая мера равна нулю. Следовательно, почти каждая точка (в
смысле меры $\mu$) нашего компакта $M$ имеет орбиту, не
пересекающуюся с такими склеенными интервалами. Мы можем построить
кусочно-непрерывное отображение на компакте $M'$ порождающем р.р.
слово $W'$, которое эквивалентно $W$ и, следовательно,  из
соображений компактности, на $M'$ есть точка, существенная эволюция
которой есть $W$. Действительно, если есть точка с существенной
эволюцией $W'$, то любое эквивалентное слово (в том числе и $W$)
также встретится.

Доказана следующая

\begin{theorem}
Множество слов, порождаемых кусочно-непрерывным преобразованием
отрезка эквивалентно множеству слов, порождаемых преобразованием
перекладывания отрезка.
\end{theorem}

Как следствие, мы получаем следующие теоремы:
\begin{theorem}
Чтобы равномерно-рекуррентное слово $W$ порождалось перекладыванием
отрезка с сохранением ориентации, необходимо и достаточно чтобы
слово обеспечивалось асимптотически правильной ориентированной
эволюцией размеченных графов Рози.
\end{theorem}

В случае меняющего ориентацию перекладывания доказана

\begin{theorem}
Чтобы равномеро-рекуррентное слово $W$ порождалось перекладыванием
отрезков, необходимо и достаточно чтобы слово обеспечивалось
асимптотически правильной эволюцией размеченных графов Рози.
\end{theorem}